\documentclass{amsart}
\usepackage[all]{xy}
\usepackage{graphics}
\usepackage{epsfig}
\usepackage{graphicx}
\usepackage{euscript,amsfonts,amssymb,amsmath,amsthm}
\usepackage[dvips]{pstricks}

 \newtheorem{thm}{Theorem}[subsection]
 \newtheorem{cor}[thm]{Corollary}

 \theoremstyle{definition}
 \newtheorem{defn}[thm]{Definition}
 \theoremstyle{remark}
 
 \numberwithin{equation}{subsection}

 \DeclareMathOperator{\Hom}{Hom}
 \DeclareMathOperator{\End}{End}
  \DeclareMathOperator{\id}{id}
 \DeclareMathOperator{\tr}{tr}

\begin{document}
\title[On Symmetrization of $6j$-symbols and Levin-Wen Hamiltonian]
{On Symmetrization of $6j$-symbols and Levin-Wen Hamiltonian}

\author{Seung-moon Hong}
\email{seungmoon.hong@utoledo.edu}
\address{Department of Mathematics\\
    The University of Toledo\\
    Toledo, OH 43606\\
    U.S.A.}

\thanks{}

\begin{abstract}

It is known that every ribbon category with unimodality allows
symmetrized $6j$-symbols with full tetrahedral symmetries while a
spherical category does not in general. We give an explicit
counterexample for this, namely the category $\mathcal{E}$. We
define the mirror conjugate symmetry of $6j$-symbols instead and
show that $6j$-symbols of any unitary spherical category can be
normalized to have this property.
 As an application, we discuss an exactly soluble model on a
honeycomb lattice. We prove that the Levin-Wen Hamiltonian is
exactly soluble and hermitian on a unitary spherical category.

\end{abstract}
\maketitle

\section{Introduction}

At first notice that in this paper $6j$-symbols with 24 full
tetrahedral symmetries are called symmetrized $6j$-symbols rather
than normalized $6j$-symbols as in \cite{T}. However they have the
same meaning.

This paper is about spherical fusion categories and their
$6j$-symbols. More specifically the $6j$-symbols in the well-known
examples are invariant under the symmetry group of the
tetrahedron. For example, this is the case if the category is
unimodal. In \cite{T}, V. Turaev showed that a ribbon category
with unimodality allows symmetrized $6j$-symbols (see chapter 6 in
\cite{T}), which implies that the state sum model on closed
3-manifold is invariant under the bistellar moves on
triangulations. However this need not be the case in general and
in section \ref{impossiblily} we show that a spherical fusion
category constructed in \cite{HH} and \cite{HRW} does not admit
$6j$-symbols with tetrahedral symmetry.

We are also interested in applications to topological field
theory. There are two standard theories associated with a
spherical category. One is the Turaev-Viro theory in \cite{TV};
Barrett and Westbury also constructed an invariant on a spherical
category in \cite{BW} and they asserted it is equivalent to the
Turaev-Viro theory. The other is the Reshetikhin-Turaev theory
based on the quantum double of a spherical category. It is an open
problem to determine if these two theories are equal. An
alternative approach is the Hamiltonian formulation of the
Turaev-Viro model in \cite{LW}. It is expected that the ground
states of this model form a modular functor and that this is
isomorphic to the Reshetikhin-Turaev theory (For this connection,
see section 6 in \cite{RSW}). In \cite{LW} it is assumed that the
$6j$-symbols have tetrahedral symmetry and it is shown that the
Hamiltonian is exactly soluble. Levin and Wen asserted in the
paper that the Hamiltonian is hermitian.

We also study unitary spherical categories and show in section
\ref{unitary-cat} that any such category admits $6j$-symbols with
so-called mirror conjugate symmetry. Then in section \ref{Hamil}
we extend the Levin-Wen Hamiltonian formulation to unitary
spherical categories with $6j$-symbols having mirror conjugate
symmetry. This Hamiltonian is also exactly soluble and hermitian.
This is a genuine extension, as the example in \cite{HH} and
\cite{HRW} is unitary.


Here are the contents of this paper. In section \ref{def}, we
recall the definition of a spherical category and make notations
for associativity and $6j$-symbols. Notice that we define two
different types of $6j$-symbols, which is unavoidable because of
the lack of symmetries as shown in section \ref{norm of E};
Barrett and Westbury did not assume symmetrized $6j$-symbols in
\cite{BW} and also defined two types of $6j$-symbols. In section
\ref{unitary-cat}, we define the mirror conjugate symmetry of
$6j$-symbols and study some other properties of a unitary
spherical category. We give normalization conditions on the
trivalent basis to obtain the mirror conjugate symmetry of
$6j$-symbols on any unitary spherical category. Section \ref{norm
of E} is devoted to symmetrized $6j$-symbols. We recall the
definition of the symmetrized $6j$-symbols and show that a
spherical category does not always allow them by giving a
counterexample, category $\mathcal{E}$ (see \cite{HH} and
\cite{HRW} for more detail on the structure of category
$\mathcal{E}$). We prove the impossibility of the symmetrization
for the category $\mathcal{E}$, and instead give a normalization
having properties studied in section \ref{unitary-cat}. In section
\ref{Hamil}, we reformulate the Levin-Wen Hamiltonian and show
that it is exactly soluble and hermitian on a unitary spherical
category. In the proof, the mirror conjugate symmetry plays an
important role along with other properties shown in section
\ref{unitary-cat}. In this sense the unitary spherical category is
good enough to define the Hamiltonian.

This paper gives only a partial result on the Hamiltonian
formulation because the study of the ground states still remains
to be done. We expect that it would be an interesting further
direction.

\section{Definitions and Notations}\label{def}

\begin{defn}
 A {\em{tensor category}} is a category $\mathcal{C}$ with a covariant functor $\otimes : \mathcal{C} \times
\mathcal{C}\rightarrow \mathcal{C}$, a natural isomorphism $\alpha
: (\mathcal{C}\otimes \mathcal{C})\otimes \mathcal{C} \rightarrow
\mathcal{C} \otimes (\mathcal{C} \otimes \mathcal{C})$, called
associativity, satisfying the pentagon axiom, a tensor unit
$\textbf{1} \in \mathcal{C}$ and natural isomorphisms $\rho :
\mathcal{C} \otimes \textbf{1}\rightarrow \mathcal{C}$ and
$\lambda : \textbf{1} \otimes \mathcal{C} \rightarrow \mathcal{C}$
satisfying the triangle axiom.

A tensor category is called {\em{strict}} if $\alpha,\rho$ and
$\lambda$ are identity.

For an algebraically closed field $k$, a tensor category is
{\em{$k$-linear}} if all $\Hom$-spaces are $k$-linear vector
spaces, composition of morphisms are $k$-bilinear, and $\otimes$
is $k$-bilinear on morphisms.

In a $k$-linear tensor category, an object $x$ in $\mathcal{C}$ is
said to be {\em{simple}} if the map $k \rightarrow \End(x)$, $c
\mapsto c\cdot \id_x$, is an isomorphism.

A $k$-linear tensor category is said to be {\em{semi-simple}} if
every object is isomorphic to a direct sum of finitely many simple
objects.

 \end{defn}

 {\em{Right rigidity}} of a strict tensor category $\mathcal{C}$ means that
for each $x \in \mathcal{C}$ there is a right dual object $x^{*}
\in \mathcal{C}$ with a morphism $b_x : \textbf{1} \rightarrow x
\otimes x^*$ and $ d_x : x^{*} \otimes x \rightarrow \textbf{1}$
such that $ (\id_x \otimes d_x)\circ (b_x \otimes \id_x)=\id_x$
and $(d_x \otimes \id_{x^*})\circ (\id_{x^*} \otimes b_x)=
\id_{x^*}$. If every object in $\mathcal{C}$ has its right dual
$*$, then we may view $*$ as a contravariant functor and $f^* \in
\Hom_{\mathcal{C}}({y^*},{x^*})$ is defined by $f^* = (d_y \otimes
\id_{x^*})  \circ (\id_{y^*} \otimes f \otimes \id_{x^*}) \circ
(\id_{y^*} \otimes b_x)$ for each $f\in
\Hom_{\mathcal{C}}({x},{y})$. The contravariant functor from
{\em{left rigidity}} is defined similarly.

A {\em{strict pivotal structure}} of a strict tensor category
$\mathcal{C}$ means that right dual is equal to left dual as a
functor. Note that if a strict category $\mathcal{C}$ has a strict
pivotal structure, then $x^{**}=x$ and $f^{**}=f$ for any object
$x$ and for any morphism $f$.

\begin{defn}
A tensor category $\mathcal{C}$ is {\em{rigid}} if it has right
and left rigidity.

A {\em{fusion category}} is a $k$-linear semi-simple rigid tensor
category with finitely many isomorphism classes of simple objects,
finite dimensional spaces of morphisms, and $\End(\textbf{1})\cong
k$.

In a strict tensor category $\mathcal{C}$ with a strict pivotal
structure, for any endomorphism $f\in \End_{\mathcal{C}}({x})$ the
{\em{right trace}} is defined by $\tr_{r}(f)=d_{x^*}\circ
(f\otimes \id_{x^*})\circ b_{x}$ and the {\em{left trace}} is
defined by $\tr_{l}(f)=d_x\circ (\id_{x^*}\otimes f)\circ
b_{x^*}$. Category $\mathcal{C}$ is {\em{spherical}} if
$\tr_{r}(f)=\tr_{l}(f)$ for all $f$, and in this case we denote it
by $\tr(f)$. In a spherical category the {\em{quantum dimension}}
$\dim (x)$ of an object $x$ is defined by $\tr(\id_x)$.

\end{defn}

In a spherical category with a strict pivotal structure, we have
the following property on trace:

\begin{align}\label{tr}
\text{$\tr(f\circ g)=\tr(g\circ f)$}
\end{align}

\noindent This equality comes from the strict pivotality as
follows:

\vspace{5mm}
\includegraphics[height=19.2mm,width=12mm]{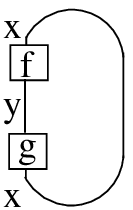}
$\begin{array}{c}=
\\ \\ \\  \end{array}$
\includegraphics[height=19.2mm,width=16.8mm]{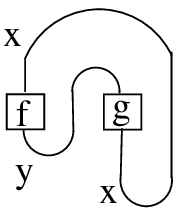}
$\begin{array}{c}=
\\ \\ \\  \end{array}$
\includegraphics[height=19.2mm,width=19.2mm]{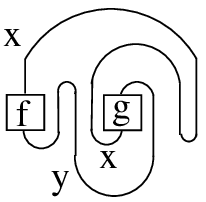}
$\begin{array}{c}=
\\ \\ \\  \end{array}$
\includegraphics[height=19.2mm,width=15.6mm]{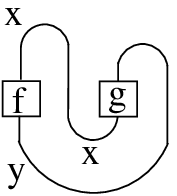}
$\begin{array}{c}=
\\ \\ \\  \end{array}$
\includegraphics[height=19.2mm,width=12mm]{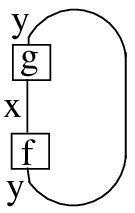}

\noindent where the second equality comes from the strict
pivotality $g=g^{**}$ and the rest are by the rigidity.

\subsection{Associativity}

In a strict spherical category $\mathcal{C}$, semi-simplicity of
the category allows every morphism to be built up from trivalent
morphisms in $\Hom_{\mathcal{C}}$-spaces of the type
$\Hom_{\mathcal{C}}(u \otimes v, w)$ or $\Hom_{\mathcal{C}}(u ,v
\otimes w)$ for simple objects $u,v$ and $w$.

Once we fix a trivalent basis for each $\Hom_{\mathcal{C}}$-space
of the type of $\Hom_{\mathcal{C}}(u \otimes v, w)$, then
associativity is represented by a matrix $F$ of the following
form:

$\begin{array}{c}F^{x}_{uvw} : \\ \\ \end{array}$
\includegraphics[height=14mm,width=16mm]{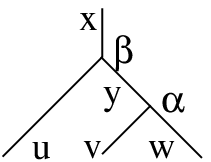}
$\begin{array}{c}=
\sum_{z,\gamma,\delta}(F^{x}_{uvw})^{y\alpha\beta}_{z\gamma\delta}
\\ \\ \end{array}$
\includegraphics[height=14mm,width=16mm]{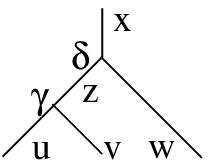}

\noindent where $\alpha,\beta,\gamma$ and $\delta$ are basis
elements in $\Hom_{\mathcal{C}}({v\otimes w},y)$, $
\Hom_{\mathcal{C}}({u\otimes y},x)$, $\Hom_{\mathcal{C}}({u\otimes
 v},z)$ and $\Hom_{\mathcal{C}}({z\otimes w},x)$, respectively.
In this convention $F^{x}_{uvw}$ is the matrix with entries
$\left(F^{x}_{uvw}\right)^{y \alpha\beta}_{z\gamma\delta}$ in the
($y \alpha\beta$)-th column and ($z\gamma\delta$)-th row.

Using a trivalent basis for each $\Hom_{\mathcal{C}}$-space of the
type of $\Hom_{\mathcal{C}}(u,  v\otimes w)$, associativity is
similarly represented by a matrix $G$ as follows:

$\begin{array}{c} G^{uvw}_{x}:\\ \\ \end{array}$
\includegraphics[height=12mm,width=20mm]{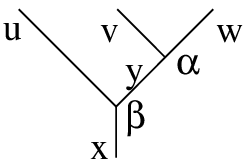}
$\begin{array}{c} = \sum_{z,\gamma,\delta}
(G^{uvw}_{x})^{y\beta\alpha}_{z\delta\gamma} \\ \\
\end{array}$
\includegraphics[height=12mm,width=20mm]{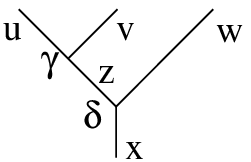}

\subsection{$6j$-symbols}

We define two different $6j$-symbols as follows :

\includegraphics[height=12mm,width=13mm]{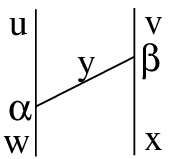}
$\begin{array}{c}= \sum_{z,\gamma,\delta}\left\{
\begin{smallmatrix}u&v&y(\alpha\beta)\\w&x&z(\gamma\delta)\end{smallmatrix} \right\}_{+} \\ \\ \end{array} \:\:$
\includegraphics[height=12mm,width=12mm]{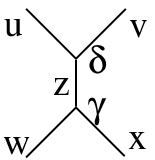}

\includegraphics[height=12mm,width=13mm]{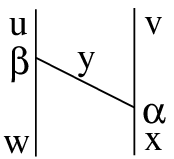}
$\begin{array}{c}= \sum_{z,\gamma,\delta}\left\{
\begin{smallmatrix}u&v&y(\alpha\beta)\\w&x&z(\gamma\delta)\end{smallmatrix} \right\}_{-} \\ \\ \end{array} \:\:$
\includegraphics[height=12mm,width=12mm]{diagrams/winfig37.eps}

\noindent In this convention, for either the $+$ or $-$ case,
$\left\{
\begin{smallmatrix}u&v&-\\w&x&-\end{smallmatrix} \right\}_{\pm}$ is the
matrix with entries $\left\{
\begin{smallmatrix}u&v&y(\alpha\beta)\\w&x&z(\gamma\delta)\end{smallmatrix}
\right\}_{\pm}$ in the ($y\alpha\beta$)-th column and
($z\gamma\delta$)-th row.

\section{Normalization of Trivalent Basis}\label{unitary-cat}

In this section, we always assume that every simple object is
self-dual and the given unitary spherical category has a strict
pivotal structure. The category is said to be unitary if all
$F$-matrices are unitary on a properly normalized trivalent basis.
In this section we show how to obtain unitary $G$-matrices and
unitary $6j$-symbols from unitary $F$-matrices by choosing
algebraic dual basis.

\subsection{Algebraic Dual Basis}

For each basis vector $\alpha \in \Hom_{\mathcal{C}}({u\otimes
v},x)$, we choose the orthogonal algebraic dual basis vector
$\bar{\alpha} \in \Hom_{\mathcal{C}}({x},u\otimes v)$ satisfying
the following condition

\begin{align}\label{theta1}
\begin{array}{c}\beta\circ \bar{\alpha} = \delta_{\alpha,\beta} \dfrac{\sqrt{uv}}{\sqrt{x}} \id_x \:\: ,\:\:\\ \\
 \end{array}
\includegraphics[height=16mm,width=12mm]{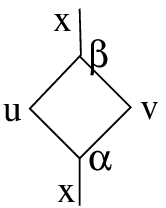}
\begin{array}{c}= \delta_{\alpha,\beta} \dfrac{\sqrt{uv}}{\sqrt{x}}
\\ \\  \end{array}
\includegraphics[height=16mm,width=4mm]{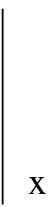}
\end{align}

\noindent Here note that $u,v$ and $x$ on the right hand side of
the equality are simple notations for the quantum dimensions of
objects $u,v$ and $x$, and that the bar on the dual basis is
omitted in the diagram whenever it is clear from the context. With
this dual basis, the following holds :

\begin{align}\label{theta2}
\left\{    \begin{array}{c}\tr(\alpha \circ
\bar{\alpha})=\sqrt{uvx}
\:\: \forall \:\:\text{basis vector}\:\: \alpha\in \Hom_{\mathcal{C}}({u\otimes v},x)  \\
\id_{u \otimes v} = \sum_{x,\alpha}\frac{\sqrt{x}}{\sqrt{uv}} \:\:
\bar{\alpha}\circ \alpha
\end{array} \right.
\end{align}

\hspace{10mm}\includegraphics[height=16mm,width=14mm]{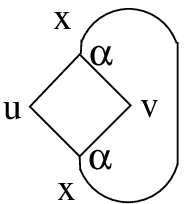}
$\begin{array}{c}= \sqrt{uvx}
\\ \\  \end{array} \:\:, \:\:$
\includegraphics[height=16mm,width=10mm]{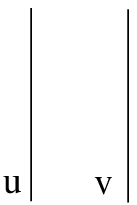}
$\begin{array}{c}= \sum_{x,\alpha}\dfrac{\sqrt{x}}{\sqrt{uv}}
\\ \\  \end{array}$
\includegraphics[height=16mm,width=12mm]{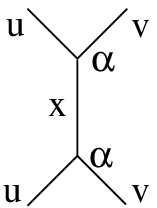}

\noindent where the summation in the second equality runs over all
simple objects $x$ and all trivalent basis morphisms $\alpha \in
\Hom_{\mathcal{C}}(u \otimes v,x)$.

Using this convention on the algebraic dual basis, $6j$-symbols
can be obtained from $F$ matrices, or vice versa, by the formulas

\begin{align}\label{6jas}
\left\{
\begin{smallmatrix}u&v&y(\alpha\beta)\\w&x&z(\gamma\delta)\end{smallmatrix} \right\}_{+}
=\dfrac{\sqrt{yz}}{\sqrt{vw}}\left(F^{z}_{uyx}\right)^{v
\beta\delta}_{w\alpha\gamma}, \left\{
\begin{smallmatrix}u&v&y(\alpha\beta)\\w&x&z(\gamma\delta)\end{smallmatrix} \right\}_{-}=
\dfrac{\sqrt{yz}}{\sqrt{ux}} \left((F^z_{wyv})^{-1}\right)^{u
\beta\delta}_{x\alpha\gamma}
\end{align}

\subsection{Mirror Conjugate Symmetry}

We may ask what happens to the transformation rules when we take
the mirror image of a given diagram about the horizontal axis.
Mirror conjugate symmetry means that in the mirror image we have
the conjugate coefficients for each transformation. It will be
shown in the next subsection that every unitary spherical category
with strict pivotality has this property. Because every
transformation of a diagram can be obtained by a sequence of
associativities or $6j-$symbols we need only to study those. In
terms of associativities, the mirror conjugate symmetry is
expressed as follows :

\vspace{3mm}

$$(G^{uvw}_{x})^{y\beta\alpha}_{z\delta\gamma} =
\overline{(F^{x}_{uvw})^{y\alpha\beta}_{z\gamma\delta}}$$

\vspace{3mm}

For the $6j$-symbols, note that the $(+)6j$-symbol and
$(-)6j$-symbol are the mirror images of each other. Thus equality
$$\left\{
\begin{smallmatrix}w&x&y(\beta\alpha)\\u&v&z(\delta\gamma)\end{smallmatrix} \right\}_{-}=\overline{\left\{
\begin{smallmatrix}u&v&y(\alpha\beta)\\w&x&z(\gamma\delta)\end{smallmatrix}
\right\}_{+}}$$ implies the mirror conjugate symmetry.

\subsection{Some Properties} We are considering a spherical category $\mathcal{C}$ in which
all associativity matrices $F$ are unitary and the algebraic dual
bases satisfy the condition (\ref{theta1}).

\begin{thm}
 We have

\begin{enumerate}

\item the category $\mathcal{C}$ has mirror conjugate symmetry,
and

\item $6j$-symbols form unitary matrices.

\end{enumerate}
\end{thm}

\begin{proof}

\begin{enumerate}
\item From the condition (\ref{theta1}) of the algebraic dual
basis,

\includegraphics[height=24mm,width=20mm]{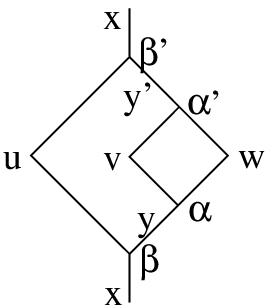}
$\begin{array}{c}=\delta_{y,y'}\delta_{\alpha,\alpha'}\delta_{\beta,\beta'}
\dfrac{\sqrt{vw}}{\sqrt{y}}\dfrac{\sqrt{uy}}{\sqrt{x}}\\ \\ \\
\end{array}$
\includegraphics[height=24mm,width=4mm]{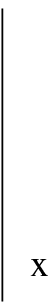}

On the other hand, using the associativities $F$ and $G$,

\vspace{10mm}
\includegraphics[height=24mm,width=20mm]{diagrams/winfig40.eps}
$\begin{array}{c}=\sum_{z,\gamma,\delta}\left(F^{x}_{uvw}\right)^{y'
\alpha'\beta'}_{z\gamma\delta}\left(G^{uvw}_{x}\right)^{y\beta\alpha}_{z\delta\gamma}\\
\\ \\
\end{array}$
\includegraphics[height=24mm,width=20mm]{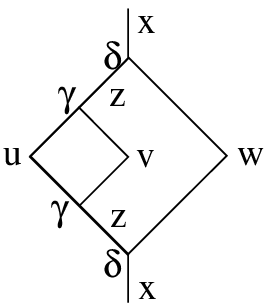}

$\begin{array}{c}=\sum_{z,\gamma,\delta}\left(F^{x}_{uvw}\right)^{y'
\alpha'\beta'}_{z\gamma\delta}\left(G^{uvw}_{x}\right)^{y\beta\alpha}_{z\delta\gamma}
\dfrac{\sqrt{uv}}{\sqrt{z}}\dfrac{\sqrt{zw}}{\sqrt{x}}\\ \\ \\
\end{array}$
\includegraphics[height=24mm,width=4mm]{diagrams/winfig41.eps}

So we have
$\delta_{y,y'}\delta_{\alpha,\alpha'}\delta_{\beta,\beta'}=\sum_{z,\gamma,\delta}\left(F^{x}_{uvw}\right)^{y'
\alpha'\beta'}_{z\gamma\delta}\left(G^{uvw}_{x}\right)^{y\beta\alpha}_{z\delta\gamma}$
while unitarity of the $F$ matrix implies
$\sum_{z,\gamma,\delta}\left(F^{x}_{uvw}\right)^{y'
\alpha'\beta'}_{z\gamma\delta}
\overline{\left(F^{x}_{uvw}\right)^{y
\alpha\beta}_{z\gamma\delta}}=\delta_{y,y'}\delta_{\alpha,\alpha'}\delta_{\beta,\beta'}$.
From the uniqueness of $F^{-1}=F^{\dag}$, we conclude
$\left(G^{uvw}_{x}\right)^{y\beta\alpha}_{z\delta\gamma}=\overline{\left(F^{x}_{uvw}\right)^{y
\alpha\beta}_{z\gamma\delta}}$.

The mirror conjugate symmetry of $6j$-symbols can be shown easily
using the formulas (\ref{6jas}),

$\left\{
\begin{smallmatrix}w&x&y(\beta\alpha)\\u&v&z(\delta\gamma)\end{smallmatrix} \right\}_{-}
=\dfrac{\sqrt{yz}}{\sqrt{vw}}
\left(\widetilde{F^{z}_{uyx}}\right)^{w
\alpha\gamma}_{v\beta\delta}=\dfrac{\sqrt{yz}}{\sqrt{vw}}
\overline{\left(F^{z}_{uyx}\right)^{v
\beta\delta}_{w\alpha\gamma}}=\overline{\left\{
\begin{smallmatrix}u&v&y(\alpha\beta)\\w&x&z(\gamma\delta)\end{smallmatrix} \right\}_{+}}.$

\vspace{3mm}

\item

For the unitarity of $6j$-symbols, we need to show
$$\sum_{z,\gamma,\delta}\left\{
\begin{smallmatrix}u&v&y(\alpha\beta)\\w&x&z(\gamma\delta)\end{smallmatrix} \right\}_{+}
\overline{\left\{
\begin{smallmatrix}u&v&y'(\alpha'\beta')\\w&x&z(\gamma\delta)\end{smallmatrix} \right\}_{+}}=
\delta_{y,y'}\delta_{\alpha,\alpha'}\delta_{\beta,\beta'}$$ and by
the mirror conjugate symmetry we need to show equivalently
$$\sum_{z,\gamma,\delta}\left\{
\begin{smallmatrix}u&v&y(\alpha\beta)\\w&x&z(\gamma\delta)\end{smallmatrix} \right\}_{+}
\left\{
\begin{smallmatrix}w&x&y'(\beta'\alpha')\\u&v&z(\delta\gamma)\end{smallmatrix} \right\}_{-}
= \delta_{y,y'}\delta_{\alpha,\alpha'}\delta_{\beta,\beta'}.$$ We
evaluate the following diagram in two different ways and then
compare them.

On one hand,

\vspace{5mm}\includegraphics[height=24mm,width=16mm]{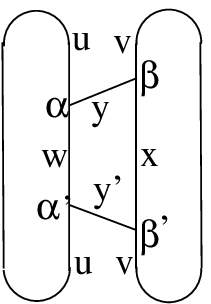}
$\begin{array}{c}\stackrel{(1)}{=}\\
\\ \\
\end{array}$
\includegraphics[height=24mm,width=16mm]{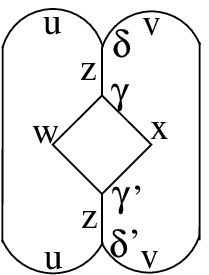}
$\begin{array}{c}\stackrel{(2)}{=}\\
\\ \\
\end{array}$
\includegraphics[height=24mm,width=18mm]{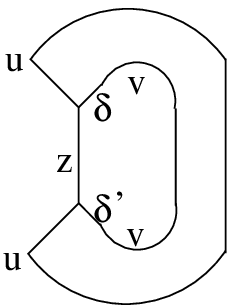}
$\begin{array}{c}\stackrel{(3)}{=}\\
\\ \\
\end{array}$
\includegraphics[height=24mm,width=14mm]{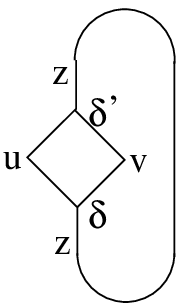}

\noindent where, for each step, the contribution is \\
$(1)=\sum_{z,\gamma,\gamma',\delta,\delta'} \left\{
\begin{smallmatrix}u&v&y(\alpha\beta)\\w&x&z(\gamma\delta)\end{smallmatrix} \right\}_{+} \left\{
\begin{smallmatrix}w&x&y'(\beta'\alpha')\\u&v&z(\delta'\gamma')\end{smallmatrix}
\right\}_{-}$ applying two $6j$-symbols,
$(2)=\delta_{\gamma,\gamma'}\dfrac{\sqrt{wx}}{\sqrt{z}}$ using
(\ref{theta1}) and sphericity. The equality $(3)$ comes from
(\ref{tr}), and the last diagram is equal to
$\delta_{\delta,\delta'} \sqrt{uvz}$ by (\ref{theta2}), thus
overall we have $\sqrt{uvwx}\sum_{z,\gamma,\delta} \left\{
\begin{smallmatrix}u&v&y(\alpha\beta)\\w&x&z(\gamma\delta)\end{smallmatrix} \right\}_{+} \left\{
\begin{smallmatrix}w&x&y'(\beta'\alpha')\\u&v&z(\delta \gamma)\end{smallmatrix}
\right\}_{-}$.

On the other hand,

\vspace{5mm}
\includegraphics[height=24mm,width=16mm]{diagrams/winfig43.eps}
$\begin{array}{c}\stackrel{(1)}{=}\\
\\ \\
\end{array}$
\includegraphics[height=24mm,width=16mm]{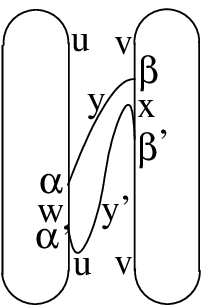}
$\begin{array}{c}\stackrel{(2)}{=}\\
\\ \\
\end{array}$
\includegraphics[height=24mm,width=20mm]{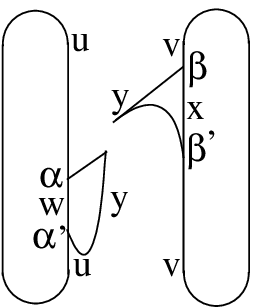}

\noindent where $(1)$ is the rigidity of $y'$,
$(2)=\delta_{y,y'}\dfrac{1}{y}$ using (\ref{theta2}), and the last
diagram is equal to
$\delta_{\alpha,\alpha'}\delta_{\beta,\beta'}\sqrt{uwy}\sqrt{yvx}$
using the sphericity, (\ref{tr}), and (\ref{theta2}), thus overall
we have
$\delta_{y,y'}\delta_{\alpha,\alpha'}\delta_{\beta,\beta'}\sqrt{uvwx}$.

\end{enumerate}
\end{proof}

\begin{cor} Let a diagram has two parts as shown below with sum over all
basis elements $\alpha$ of $\Hom_{\mathcal{C}}(u\otimes x,{v})$
(or $\beta$ of $\Hom_{\mathcal{C}}(x\otimes v,{u})$). Then we can
transform the two parts in the diagram simultaneously with sum
over all basis $\gamma$ of $\Hom_{\mathcal{C}}({u\otimes v},x)$ as
follows:

$$
\begin{array}{c}\sum_{\alpha} \\ \\ \\ \end{array}
\text{\includegraphics[height=18mm,width=36mm]{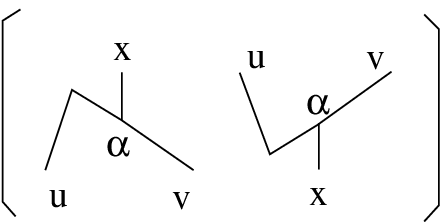}}
\begin{array}{c}=\sum_{\gamma} \\ \\ \\ \end{array}
\text{\includegraphics[height=18mm,width=36mm]{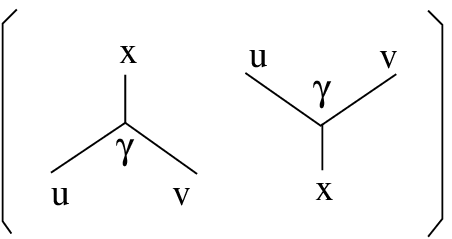}}
\begin{array}{c}=\sum_{\beta} \\ \\ \\ \end{array}
\text{\includegraphics[height=18mm,width=36mm]{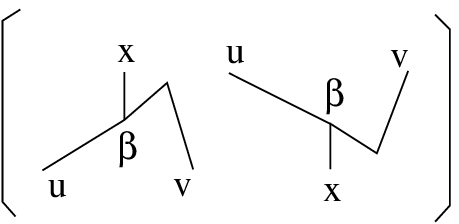}}$$

\end{cor}

\begin{proof}
This is a direct consequence of the mirror conjugate symmetry and
unitarity of $6j$-symbols. For the first equality, the coefficient
of the transformation is

$\sum_{\alpha,\gamma,\delta} \left\{
\begin{smallmatrix}\textbf{1}&x&u(\alpha)\\u&v&x(\gamma)\end{smallmatrix} \right\}_{-}
\left\{
\begin{smallmatrix}u&v&u(\alpha)\\\textbf{1}&x&x(\delta)\end{smallmatrix} \right\}_{+}$
$=\sum_{\gamma,\delta} \left(\sum_{\alpha}\overline{\left\{
\begin{smallmatrix}u&v&u(\alpha)\\\textbf{1}&x&x(\gamma)\end{smallmatrix}
\right\}_{+}} \left\{
\begin{smallmatrix}u&v&u(\alpha)\\\textbf{1}&x&x(\delta)\end{smallmatrix}
\right\}_{+}\right)=\delta_{\gamma,\delta}\sum_{\gamma} 1 $.

The proof for the second one is similar.
\end{proof}

In particular, this corollary implies that $\id_{u\otimes v}$ can
be expressed in a different way from the second equality of
(\ref{theta2}), as follows:

\begin{align}\label{theta3}
\includegraphics[height=16mm,width=9mm]{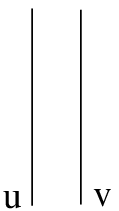}
\begin{array}{c}= \sum_{x,\gamma}\dfrac{\sqrt{x}}{\sqrt{uv}}
\\ \\ \\ \end{array}
\includegraphics[height=17mm,width=13mm]{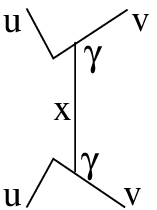}
\begin{array}{c} \:\:\: \text{or} \:\:\: \sum_{x,\gamma}\dfrac{\sqrt{x}}{\sqrt{uv}}
\\ \\ \\ \end{array}
\includegraphics[height=17mm,width=13mm]{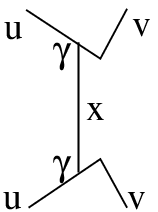}
\end{align}

Another consequence of the unitarity and mirror conjugate symmetry
is the following corollary which plays an important role when we
discuss the Levin-Wen Hamiltonian in the next section.

\begin{cor}\label{can} The following equality holds:

\vspace{5mm} \hspace{30mm} $\begin{array}{c}
\frac{\sqrt{g'}}{\sqrt{g}} \sum_{s,\eta} \sqrt{s}
\\ \\ \end{array}$
\includegraphics[height=10mm,width=24mm]{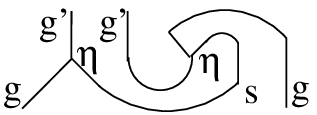}
 $\begin{array}{c} =g'
\\ \\ \end{array}$
\includegraphics[height=8mm,width=14mm]{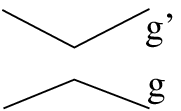}

\noindent where the summation runs over all simple objects $s$ and
all trivalent basis morphisms $\eta \in
\Hom_{\mathcal{C}}({g\otimes s},g')$.

\end{cor}

\begin{proof}
Since $**=\id$ on morphisms, we may replace the right half of the
given diagram by its double dual,

\vspace{3mm}\hspace{10mm}\includegraphics[height=16mm,width=24mm]{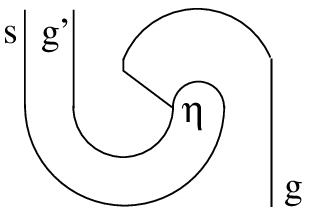}
$\begin{array}{c}=\\ \\ \\ \end{array} \:\:$
\includegraphics[height=16mm,width=19mm]{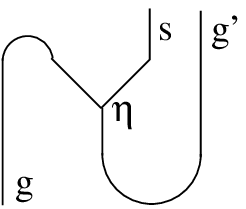}

\noindent then we set

$\begin{array}{c}\frac{\sqrt{g'}}{\sqrt{g}} \sum_{s,\eta} \sqrt{s}
\\ \\ \\  \end{array}$
\includegraphics[height=18mm,width=25mm]{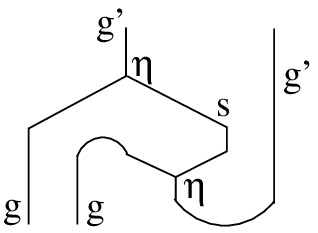}
$\begin{array}{c}= \sum_{t,\gamma,\delta} C_{t(\gamma,\delta)} \\ \\
\\
\end{array}$
\includegraphics[height=16mm,width=13mm]{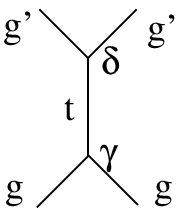}

\noindent and take the trace after composing $\gamma \in
\Hom_{\mathcal{C}}({t},g \otimes g)$ and $\delta \in
\Hom_{\mathcal{C}}({g'\otimes g'},t)$ for each $t,\gamma$ and
$\delta$.

\vspace{5mm}$\begin{array}{c}\frac{\sqrt{g'}}{\sqrt{g}}
\sum_{s,\eta} \sqrt{s}
\\ \\ \\ \\  \end{array}$
\includegraphics[height=30mm,width=28mm]{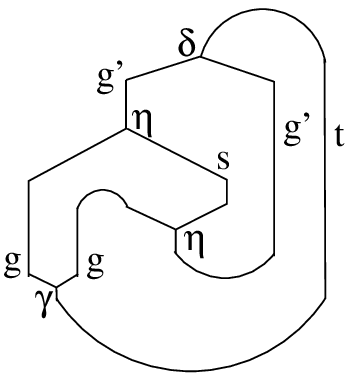}
$\begin{array}{c}= \:\:\: C_{t(\gamma,\delta)} \\ \\ \\ \\
\end{array}$
\includegraphics[height=30mm,width=21mm]{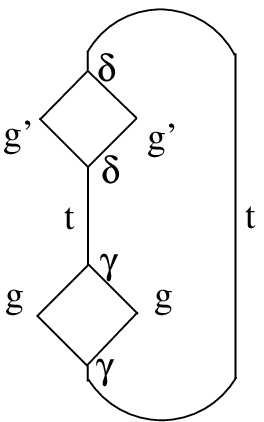}

Now we claim that $C_{t(\gamma,\delta)}=g'$ only for the case
$t=\textbf{1}$ and is $0$ otherwise. This can be shown by
deforming the diagram on the left hand side as follows (dotted box
indicates the place deformed for each step) :

\vspace{5mm}\includegraphics[height=30mm,width=28mm]{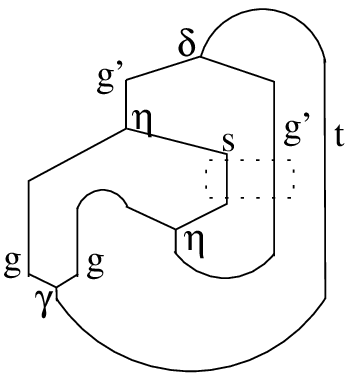}
$\begin{array}{c}\stackrel{(1)}{=}
\\ \\ \\ \\  \end{array}$
\includegraphics[height=30mm,width=28mm]{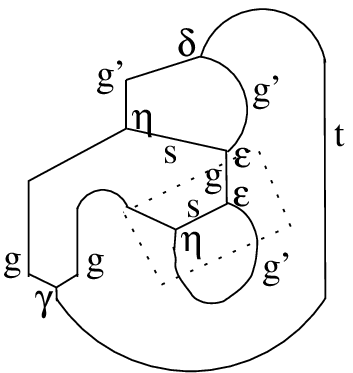}
$\begin{array}{c}\stackrel{(2)}{=}
\\ \\ \\ \\  \end{array}$
\includegraphics[height=30mm,width=28mm]{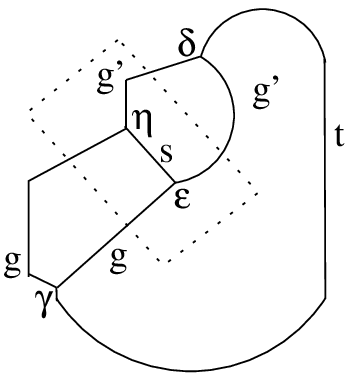}
$\begin{array}{c}\stackrel{(3)}{=} \\ \\ \\ \\  \end{array}$
\includegraphics[height=30mm,width=21mm]{diagrams/winfig62.eps}

\noindent where for each step, the contribution is
$(1)=\sum_{\varepsilon}\dfrac{\sqrt{g}}{\sqrt{sg'}}$ using
(\ref{theta2}), $(2)=g' \left\{
\begin{smallmatrix}g&g&s(\eta \varepsilon)\\g'&g'&\textbf{1} \end{smallmatrix}
\right\}_{+}$, and $(3)=\left\{
\begin{smallmatrix}g'&g'&s(\varepsilon \eta)\\g&g&t(\gamma\delta) \end{smallmatrix} \right\}_{-}$.
Thus with the initial coefficient $\frac{\sqrt{g'}}{\sqrt{g}}
\sum_{s,\eta} \sqrt{s}$ we have

$C_{t(\gamma,\delta)}=g' \sum_{s,\eta,\varepsilon}\left\{
\begin{smallmatrix}g&g&s(\eta\varepsilon)\\g'&g'&\textbf{1} \end{smallmatrix} \right\}_{+}
\left\{
\begin{smallmatrix}g'&g'&s(\varepsilon\eta)\\g&g&t(\gamma\delta) \end{smallmatrix} \right\}_{-}$
$=g' \sum_{s,\eta,\varepsilon}\left\{
\begin{smallmatrix}g&g&s(\eta\varepsilon)\\g'&g'&\textbf{1} \end{smallmatrix} \right\}_{+}
  \overline{\left\{
\begin{smallmatrix}g&g&s(\eta\varepsilon)\\g'&g'&t(\delta\gamma) \end{smallmatrix} \right\}_{+}}$
$=\delta_{t,\textbf{1}} \cdot g'$

\noindent where the second equality comes from the mirror
conjugate symmetry and the last one from the unitarity of
$6j$-symbols.

\end{proof}

\section{Non-symmetrized $6j$-symbols}\label{norm of
E}

  In
\cite{T}, V. Turaev showed that a ribbon category with so-called
unimodality allows symmetrized $6j$-symbols (see chapter 6 in
\cite{T}). Braiding in a ribbon category plays an important role
to prove this. A question now is whether or not this is possible
for a spherical category which does not allow braiding in general.
In this section, we give a counterexample to this question, which
is the category $\mathcal{E}$. The category $\mathcal{E}$ in
\cite{HH} is an example of a unitary spherical category with
strict pivotal structure (all $F$-matrices in \cite{HH} are given
again in the appendix). In the following we show that the category
$\mathcal{E}$ does not allow the symmetrized $6j$-symbols.
Furthermore, we normalize the trivalent basis morphisms in each
$\Hom_{\mathcal{E}}$ space to have the properties in section
\ref{unitary-cat} including unitary $F$-matrices, while
$F$-matrices in \cite{HH} are not yet unitary.

Notice that here we use the left multiplication convention instead
of the right multiplication convention used in \cite{HH}. So we
need to consider the transpose of each $F$-matrix in \cite{HH}.

\subsection{Symmetrized $6j$-symbols}

If the $(+)6j$-symbol is equal to the $(-)6j$-symbol, we define a
new $6j$-symbol by

$$\left\{\begin{smallmatrix}u&v&y(\alpha\beta)\\w&x&z(\gamma\delta)\end{smallmatrix}
\right\}=
\left\{\begin{smallmatrix}u&v&y(\alpha\beta)\\w&x&z(\gamma\delta)\end{smallmatrix}
\right\}_{+} =
\left\{\begin{smallmatrix}u&v&y(\beta\alpha)\\w&x&z(\gamma\delta)\end{smallmatrix}
\right\}_{-} $$

For the symmetrized $6j$-symbol we reqire 24 tetrahedral
symmetries generated by the following:

$$\left\{\begin{matrix}u&v&y\\w&x&z\end{matrix} \right\}=
\left\{\begin{matrix}v&u&y\\x&w&z\end{matrix} \right\}=
\left\{\begin{matrix}x&w&y\\v&u&z\end{matrix} \right\}=
\left\{\begin{matrix}y&x&w\\u&z&v\end{matrix} \right\}
\dfrac{\sqrt{yz}}{\sqrt{vw}}$$

\noindent where the above is a simple expression for the case of
1-dimensional $\Hom_{\mathcal{C}}$-spaces. For multi dimensional
$\Hom_{\mathcal{C}}$-space case, we need to specify trivalent
basis as before.

For a unimodal ribbon category, it is always possible for us to
have symmetrized $6j$-symbols on an appropriate basis by \cite{T}.

\subsection{Impossibility of Symmetrization for the Category
$\mathcal{E}$}\label{impossiblily}

Symmetrized 6j symbols have property $\left\{
\begin{smallmatrix}a&d&e\\b&c&-\end{smallmatrix}
\right\}_{+}=\left\{
\begin{smallmatrix}a&d&e\\b&c&-\end{smallmatrix}
\right\}_{-}$. Otherwise the horizontal line in the middle does
not have any meaning.

\vspace{3mm}\hspace{30mm}
 \psset{unit=3mm}
\begin{pspicture}(0,0)(2,3)
\psline(0,0)(0,3) \psline(2,0)(2,3)\psline(0,1)(2,2)
\put(-0.9,0){$b$}\put(-0.9,2.5){$a$}\put(2.2,0){$c$}\put(2.2,2.5){$d$}\put(0.7,1.8){$e$}
\end{pspicture}
$\begin{array}{c} \:\:\:\:=\:\:\:\: \\ \\
\\\end{array}$
\begin{pspicture}(0,0)(2,3)
\psline(0,0)(0,3) \psline(2,0)(2,3)\psline(0,2)(2,1)
\put(-0.9,0){$b$}\put(-0.9,2.5){$a$}\put(2.2,0){$c$}\put(2.2,2.5){$d$}\put(0.7,1.8){$e$}
\end{pspicture}

\vspace{3mm} The category $\mathcal{E}$, however, does not allow
it. In other words, no matter how one normalizes the trivalent
basis, the equality can not be obtained. Supposing it, we easily
get a contradiction as follows:

Suppose we have new basis for each $\Hom_{\mathcal{E}}$ space with
such a property. Then we express each side of the above equality
as a linear combination of the old basis elements and then compare
them.

Let $w^{\textbf{1}}_{xx}$, $w^{xx}_{\textbf{1}}$, $ \{w_{1}$, $
w_{2} \}$, $ \{w^{1}$, $ w^{2} \}$ be new bases of
$\Hom_{\mathcal{E}}({x\otimes x},{\textbf{1}})$, $
\Hom_{\mathcal{E}}({\textbf{1}},{x\otimes x})$, $
\Hom_{\mathcal{E}}({x\otimes x},x)$, $
\Hom_{\mathcal{E}}({x},x\otimes x)$, respectively, and let
$w^{\textbf{1}}_{xx}=f\cdot v^{\textbf{1}}_{xx}$, $
w^{xx}_{\textbf{1}}=f'\cdot v^{xx}_{\textbf{1}}$, $w_{1}=k\cdot
v_{1}+l\cdot v_{2}$, $ w_{2} = m\cdot v_{1}+n\cdot v_{2}$ ,
$w^{1}=k'\cdot v^{1}+l'\cdot v^{2}$, $ w^{2} = m'\cdot
v^{1}+n'\cdot v^{2}$ for some nonzero $f, f'$ and invertible
$\left[
\begin{smallmatrix} k&m \\ l&n \end{smallmatrix} \right], \left[
\begin{smallmatrix} k'&m' \\ l'&n' \end{smallmatrix} \right]$
where the basis elements denoted by $v$ are the old basis used in
\cite{HH}. In the following diagram, the thickened graphs denote
these new basis elements. Note that we are using the same
convention for diagrams as in \cite{HH}, that is, the solid line
indicates the object $x$, and dotted line indicates the object
$y$.

\vspace{3mm} \psset{unit=3mm,linewidth=1pt}
\begin{pspicture}(0,0)(2,1)
\psline(0,0)(1,1) \psline(1,1)(2,0)
\end{pspicture}
 = $f$
\psset{unit=3mm,linewidth=0.5pt}
\begin{pspicture}(0,0)(2,1)
\psline(0,0)(1,1) \psline(1,1)(2,0)
\end{pspicture}\hspace{5mm},\hspace{5mm}
\psset{unit=3mm,linewidth=1pt}
\begin{pspicture}(0,0)(2,1)
\psline(0,1)(1,0) \psline(1,0)(2,1)
\end{pspicture}
 = $f'$
\psset{unit=3mm,linewidth=0.5pt}
\begin{pspicture}(0,0)(2,1)
\psline(0,1)(1,0) \psline(1,0)(2,1)
\end{pspicture}

\vspace{3mm}
 \psset{unit=3mm,linewidth=1pt}
\begin{pspicture}(0,0)(2,2)
\psline(1,2)(1,1)\psline(0,0)(1,1) \psline(1,1)(2,0)
\put(0.65,0.65){$\bullet$}
\end{pspicture}
= $k$ \psset{unit=3mm,linewidth=0.5pt}
\begin{pspicture}(0,0)(2,2)
\psline(1,2)(1,1)\psline(0,0)(1,1) \psline(1,1)(2,0)
\put(0.65,0.65){$\bullet$}
\end{pspicture}
$+ l$ \psset{unit=3mm,linewidth=0.5pt}
\begin{pspicture}(0,0)(2,2)
\psline(1,2)(1,1)\psline(0,0)(1,1) \psline(1,1)(2,0)
\put(0.65,0.65){$\circ$}
\end{pspicture}\hspace{5mm},\hspace{5mm}
\psset{unit=3mm,linewidth=1pt}
\begin{pspicture}(0,0)(2,2)
\psline(1,2)(1,1)\psline(0,0)(1,1) \psline(1,1)(2,0)
\put(0.65,0.65){$\circ$}
\end{pspicture}
= $m$ \psset{unit=3mm,linewidth=0.5pt}
\begin{pspicture}(0,0)(2,2)
\psline(1,2)(1,1)\psline(0,0)(1,1) \psline(1,1)(2,0)
\put(0.65,0.65){$\bullet$}
\end{pspicture}
$+ n$ \psset{unit=3mm,linewidth=0.5pt}
\begin{pspicture}(0,0)(2,2)
\psline(1,2)(1,1)\psline(0,0)(1,1) \psline(1,1)(2,0)
\put(0.65,0.65){$\circ$}
\end{pspicture}

\vspace{3mm}
 \psset{unit=3mm,linewidth=1pt}
\begin{pspicture}(0,0)(2,2)
\psline(1,0)(1,1)\psline(0,2)(1,1) \psline(1,1)(2,2)
\put(0.65,0.65){$\bullet$}
\end{pspicture}
= $k'$ \psset{unit=3mm,linewidth=0.5pt}
\begin{pspicture}(0,0)(2,2)
\psline(1,0)(1,1)\psline(0,2)(1,1) \psline(1,1)(2,2)
\put(0.65,0.65){$\bullet$}
\end{pspicture}
$+ l'$ \psset{unit=3mm,linewidth=0.5pt}
\begin{pspicture}(0,0)(2,2)
\psline(1,0)(1,1)\psline(0,2)(1,1) \psline(1,1)(2,2)
\put(0.65,0.65){$\circ$}
\end{pspicture}\hspace{5mm},\hspace{5mm}
\psset{unit=3mm,linewidth=1pt}
\begin{pspicture}(0,0)(2,2)
\psline(1,0)(1,1)\psline(0,2)(1,1) \psline(1,1)(2,2)
\put(0.65,0.65){$\circ$}
\end{pspicture}
= $m'$ \psset{unit=3mm,linewidth=0.5pt}
\begin{pspicture}(0,0)(2,2)
\psline(1,0)(1,1)\psline(0,2)(1,1) \psline(1,1)(2,2)
\put(0.65,0.65){$\bullet$}
\end{pspicture}
$+ n'$ \psset{unit=3mm,linewidth=0.5pt}
\begin{pspicture}(0,0)(2,2)
\psline(1,0)(1,1)\psline(0,2)(1,1) \psline(1,1)(2,2)
\put(0.65,0.65){$\circ$}
\end{pspicture}

\vspace{3mm}

Now consider the case that $a = \mathbf{1}$ and everything else is
$x$ in the above symmetry.

\vspace{3mm} $\big( \:\: $\psset{unit=3mm,linewidth=1pt}
\begin{pspicture}(0,0)(2,2)
\psline(0,0)(0,0.7) \psline(2,0)(2,2)\psline(0,0.7)(2,1.3)
\put(1.65,0.85){$\bullet$}
\end{pspicture}  ,
\begin{pspicture}(0,0)(2,2)
\psline(0,0)(0,0.7) \psline(2,0)(2,2)\psline(0,0.7)(2,1.3)
\put(1.65,0.85){$\circ$}
\end{pspicture}$\:\: \big)$
$\;\;\;\;$=$\;\;\;\;$ $\big( \:\: $\psset{unit=3mm,linewidth=1pt}
\begin{pspicture}(0,0)(2,2)
\psline(0,0)(0,1.3) \psline(2,0)(2,2)\psline(0,1.3)(2,0.7)
\put(1.65,0.4){$\bullet$}
\end{pspicture}  ,
\begin{pspicture}(0,0)(2,2)
\psline(0,0)(0,1.3) \psline(2,0)(2,2)\psline(0,1.3)(2,0.7)
\put(1.65,0.4){$\circ$}
\end{pspicture}$\:\: \big)$

\vspace{3mm}

If we express the left hand side in terms of the old basis $v_1,
v_2$, then we have a matrix expression $\left[
\begin{smallmatrix} k&m \\ l&n \end{smallmatrix} \right]$. For the
right hand side, we need to use associativity $(F^{x}_{xxx})^{-1}$
to get $\frac{f}{2\dim(x)} \left[
\begin{smallmatrix} 1&1 \\ 1&-1 \end{smallmatrix} \right] \left[
\begin{smallmatrix} k'&m' \\ l'&n' \end{smallmatrix} \right]$ on the
same old basis. Thus the symmetry means $\left[
\begin{smallmatrix} k&m \\ l&n \end{smallmatrix} \right]= \frac{f}{2\dim(x)}
\left[
\begin{smallmatrix} 1&1 \\ 1&-1 \end{smallmatrix} \right]\left[
\begin{smallmatrix} k'&m' \\ l'&n' \end{smallmatrix} \right]$,
or equivalently $\left[
\begin{smallmatrix} k'&m' \\ l'&n' \end{smallmatrix} \right]= \frac{\dim(x)}{f}
\left[
\begin{smallmatrix} 1&1 \\ 1&-1 \end{smallmatrix} \right]\left[
\begin{smallmatrix} k&m \\ l&n \end{smallmatrix} \right]$.

Next, consider the case that $c=\mathbf{1}$ and everything else is
$x$.

\vspace{3mm} $\big( \:\: $\psset{unit=3mm,linewidth=1pt}
\begin{pspicture}(0,0)(2,2)
\psline(0,0)(0,2) \psline(2,1.3)(2,2)\psline(0,0.7)(2,1.3)
\put(-0.3,0.4){$\bullet$}
\end{pspicture}  ,
\begin{pspicture}(0,0)(2,2)
\psline(0,0)(0,2) \psline(2,1.3)(2,2)\psline(0,0.7)(2,1.3)
\put(-0.3,0.4){$\circ$}
\end{pspicture}$\:\: \big)$
$\;\;\;\;$=$\;\;\;\;$ $\big( \:\: $\psset{unit=3mm,linewidth=1pt}
\begin{pspicture}(0,0)(2,2)
\psline(0,0)(0,2) \psline(2,0.7)(2,2)\psline(0,1.3)(2,0.7)
\put(-0.3,0.85){$\bullet$}
\end{pspicture}  ,
\begin{pspicture}(0,0)(2,2)
\psline(0,0)(0,2) \psline(2,0.7)(2,2)\psline(0,1.3)(2,0.7)
\put(-0.3,0.85){$\circ$}
\end{pspicture}$\:\: \big)$

\vspace{3mm}

From this symmetry, we have $\left[
\begin{smallmatrix} k'&m' \\ l'&n' \end{smallmatrix} \right]=
f'e^{5\pi i/6}
  \left[
\begin{smallmatrix} 1&-i \\ -i&1 \end{smallmatrix}\right] \left[
\begin{smallmatrix} k&m \\ l&n \end{smallmatrix} \right]$.
Comparing the two matrix equalities, we have $\dfrac{\dim(x)}{f}
\left[\begin{smallmatrix} 1&1 \\ 1&-1
\end{smallmatrix} \right]=f'e^{5\pi i/6} \left[\begin{smallmatrix} 1&-i \\ -i&1 \end{smallmatrix} \right]$,
and it is easy to see the impossibility.

\subsection{Normalization for the Category
$\mathcal{E}$}

Fix $d=\dim(x)=1+\sqrt{3}$ and $v=\sqrt{d}$, and then we define
new basis elements(thick ones) as follows:

\vspace{3mm}
 \psset{unit=3mm,linewidth=1pt}
\begin{pspicture}(0,0)(2,2)
\psline(1,2)(1,1)\psline(0,0)(1,1) \psline(1,1)(2,0)
\put(1.3,1.1){$_j$}
\end{pspicture}
= $\sqrt{v}$ \psset{unit=3mm,linewidth=0.5pt}
\begin{pspicture}(0,0)(2,2)
\psline(1,2)(1,1)\psline(0,0)(1,1) \psline(1,1)(2,0)
\put(1.3,1.1){$_j$}
\end{pspicture}$\:\:$,$\:\:$
\psset{unit=3mm,linewidth=1pt}
\begin{pspicture}(0,0)(2,2)
\psline(1,0)(1,1)\psline(0,2)(1,1) \psline(1,1)(2,2)
\put(1.3,0.8){$_j$}
\end{pspicture}
=$\sqrt{v}$ \psset{unit=3mm,linewidth=0.5pt}
\begin{pspicture}(0,0)(2,2)
\psline(1,0)(1,1)\psline(0,2)(1,1) \psline(1,1)(2,2)
\put(1.3,0.8){$_j$}
\end{pspicture}$\:\:$ for each $j \in \{1, 2 \}$

\vspace{3mm}
 \psset{unit=3mm,linewidth=1pt}
\begin{pspicture}(0,0)(2,1)
\psline(0,0)(1,1) \psline(1,1)(2,0)
\end{pspicture}
 = $\sqrt{2}d$
\psset{unit=3mm,linewidth=0.5pt}
\begin{pspicture}(0,0)(2,1)
\psline(0,0)(1,1) \psline(1,1)(2,0)
\end{pspicture}$\:\:$,$\:\:$
\psset{unit=3mm,linewidth=1pt}
\begin{pspicture}(0,0)(2,1)
\psline(0,1)(1,0) \psline(1,0)(2,1)
\end{pspicture}
 = $\frac{1}{\sqrt{2}}$
\psset{unit=3mm,linewidth=0.5pt}
\begin{pspicture}(0,0)(2,1)
\psline(0,1)(1,0) \psline(1,0)(2,1)
\end{pspicture}$\:\:$,$\:\:$
\psset{unit=3mm,linewidth=1pt}
\begin{pspicture}(0,0)(2,2)
\psline[linestyle=dotted](1,2)(1,1) \psline(0,0)(1,1)
\psline(1,1)(2,0)
\end{pspicture}
= $\sqrt{2}d$ \psset{unit=3mm,linewidth=0.5pt}
\begin{pspicture}(0,0)(2,2)
\psline[linestyle=dotted](1,2)(1,1)\psline(0,0)(1,1)
\psline(1,1)(2,0)
\end{pspicture}$\:\:$,$\:\:$
\psset{unit=3mm,linewidth=1pt}
\begin{pspicture}(0,0)(2,2)
\psline[linestyle=dotted](1,0)(1,1)\psline(0,2)(1,1)
\psline(1,1)(2,2)
\end{pspicture}
=$\frac{1}{\sqrt{2}}$ \psset{unit=3mm,linewidth=0.5pt}
\begin{pspicture}(0,0)(2,2)
\psline[linestyle=dotted](1,0)(1,1)\psline(0,2)(1,1)
\psline(1,1)(2,2)
\end{pspicture}

\vspace{3mm} and keep \psset{unit=3mm,linewidth=0.5pt}
\begin{pspicture}(0,0)(2,2)
\psline(1,2)(1,1)\psline[linestyle=dotted](0,0)(1,1)
\psline(1,1)(2,0)
\end{pspicture}$\:\:$,$\:\:$
\begin{pspicture}(0,0)(2,2)
\psline(1,2)(1,1)\psline(0,0)(1,1)
\psline[linestyle=dotted](1,1)(2,0)
\end{pspicture}$\:\:$,$\:\:$
\begin{pspicture}(0,0)(2,2)
\psline(1,0)(1,1)\psline[linestyle=dotted](0,2)(1,1)
\psline(1,1)(2,2)
\end{pspicture}
$\:\:$,$\:\:$
\begin{pspicture}(0,0)(2,2)
\psline(1,0)(1,1)\psline(0,2)(1,1)
\psline[linestyle=dotted](1,1)(2,2)
\end{pspicture}$\:\:$,$\:\:$
\begin{pspicture}(0,0)(2,1)
\psline[linestyle=dotted](0,0)(1,1)
\psline[linestyle=dotted](1,1)(2,0)
\end{pspicture}$\:\:$,$\:\:$
\begin{pspicture}(0,0)(2,1)
\psline[linestyle=dotted](0,1)(1,0)
\psline[linestyle=dotted](1,0)(2,1)
\end{pspicture}
as before. \vspace{3mm}

It is a straightforward computation to see that these basis
elements allow all properties in section \ref{unitary-cat},
especially the property \ref{theta1}. We list the associativity
matrices $F$ obtained by this normalization below, all of which
are unitary. The $6j$-symbols can be computed by the formulas
\ref{6jas}

\vspace{3mm}

$F^{y}_{yyy}=$ $F^{x}_{xyy}=$ $F^{x}_{yyx}=$
$F^{\textbf{1}}_{xyx}=$ $F^{\textbf{1}}_{xxy}=$ $F^{y}_{xxy}=$
$F^{\textbf{1}}_{yxx}=$ $F^{y}_{yxx}=1$,

$F^{y}_{xyx}=$ $F^{x}_{yxy}= -1$,

$F^{x}_{xyx}= \left[\begin{smallmatrix}1&0\\0&-1\end{smallmatrix}
\right]$, $F^{x}_{xxy}=
\left[\begin{smallmatrix}0&-i\\i&0\end{smallmatrix} \right]$,
$F^{x}_{yxx}= \left[\begin{smallmatrix}0&1\\1&0\end{smallmatrix}
\right]$,

$F^{\textbf{1}}_{xxx}= \frac{1}{\sqrt{2}} e^{7 \pi i
 /12} \left[\begin{smallmatrix}1&1\\i&-i\end{smallmatrix}
\right]$, $F^{y}_{xxx}= \frac{1}{\sqrt{2}} e^{7 \pi i
 /12} \left[\begin{smallmatrix}i&-i\\1&1\end{smallmatrix}
\right]$,

$F^{x}_{xxx}=$ $\left[ \begin{smallmatrix}
\frac{1}{d}&\frac{1}{d}&\frac{1}{\sqrt{2}v}&\frac{1}{\sqrt{2}v}&\frac{1}{\sqrt{2}v}&\frac{-1}{\sqrt{2}v}\\
\frac{1}{d}&\frac{-1}{d}&\frac{1}{\sqrt{2}v}&\frac{1}{\sqrt{2}v}&\frac{-1}{\sqrt{2}v}&\frac{1}{\sqrt{2}v}\\
\frac{1}{\sqrt{2}v} e^{\frac{-5\pi i}{6}}&\frac{1}{\sqrt{2}v}
e^{\frac{-5\pi i}{6}}&\frac{1}{\sqrt{2}d}
e^{\frac{-5\pi i}{12}}&\frac{1}{2} e^{\frac{\pi i}{3}}&\frac{1}{\sqrt{2}d} e^{\frac{-5\pi i}{12}}&\frac{1}{2} e^{\frac{-2\pi i}{3}}\\
\frac{1}{\sqrt{2}v} e^{\frac{-\pi i}{3}}&\frac{1}{\sqrt{2}v}
e^{\frac{-\pi i}{3}}&\frac{1}{2} e^{\frac{5\pi
i}{6}}&\frac{1}{\sqrt{2}d} e^{\frac{\pi i}{12}}&\frac{1}{2}
e^{\frac{5\pi i}{6}}&\frac{1}{\sqrt{2}d} e^{\frac{-11\pi
i}{12}}\\
\frac{1}{\sqrt{2}v} e^{\frac{-\pi i}{3}}&\frac{1}{\sqrt{2}v}
e^{\frac{2\pi i}{3}}&\frac{1}{\sqrt{2}d} e^{\frac{\pi
i}{12}}&\frac{1}{2} e^{\frac{5\pi i}{6}}&\frac{1}{\sqrt{2}d}
e^{\frac{-11\pi
i}{12}}&\frac{1}{2} e^{\frac{5\pi i}{6}}\\
\frac{1}{\sqrt{2}v} e^{\frac{-5\pi i}{6}}&\frac{1}{\sqrt{2}v}
e^{\frac{\pi i}{6}}&\frac{1}{2} e^{\frac{\pi
i}{3}}&\frac{1}{\sqrt{2}d} e^{\frac{-5\pi i}{12}}&\frac{1}{2}
e^{\frac{-2\pi i}{3}}&\frac{1}{\sqrt{2}d} e^{\frac{-5\pi i}{12}}
\end{smallmatrix}
\right]$

\section{Exactly Soluble Lattice Models}\label{Hamil}

The exactly soluble Hamiltonian on a honeycomb lattice model is
studied in \cite{LW} and in chapter 11 of \cite{We}, in which they
assume the symmetrized $6j$-symbols. They claim that the
symmetrized $6j$-symbols with a so-called unitary condition imply
that the Hamiltonian is exactly soluble and hermitian. In this
section we study the same model but on a unitary spherical
category.

We prove the following theorem in this section:

\begin{thm}
On a unitary spherical category with strict pivotal structure, the
Levin-Wen Hamiltonian $H$ on the honeycomb lattice model has the
following properties:
\begin{enumerate}
\item $B_P$'s and $E_I$'s commute with each other and hence $H$ is
exactly soluble. \item $B_P$'s are projectors. \item $H$ is
hermitian.
\end{enumerate}
\end{thm}

Commutativity is clear. In particular the commutativity of
$B_{P_1}$ and $B_{P_2}$ is an easy conclusion from the topological
consideration. We prove the rest of the theorem in section
\ref{proj} and \ref{herm}.

\subsection{Definitions}

The Levin-Wen model is the usual spin model on the honeycomb
lattice with edges decorated by simple objects of a category
$\mathcal{C}$. In the case that any fusion coefficient
$N^{c}_{a,b}=\dim \Hom_{\mathcal{C}}(a\otimes b,c)$ is bigger than
1, we need to distinguish the corresponding vertex by decorating
it with basis morphisms of the $\Hom_{\mathcal{C}}$-space, so the
Hilbert space is $\otimes_{\text{edges}}\mathbb{C}^{n}
\otimes_{\text{vertices}} \mathbb{C}^{N^{c}_{a,b}}$ where $n$ is
the rank of the given category and three edges colored by simple
objects $a, b$, and $c$ meet at a vertex.

The exactly soluble Hamiltonian on the honeycomb lattice model is
defined by
$$H=\sum_{I}\left(1-E_I \right) +\sum_{P}\left(1-B_P \right),\:\:
B_P=\sum_s \frac{s}{D^2}B^{s}_{P} $$

\noindent $B^{s}_{P}$ is an operator acting on 12 links and 6
trivalent vertices around the hexagon $P$ by introducing an extra
loop labeled by simple object $s$ in $P$. $E_I$ acts on the vertex
$I$ such that it is the identity if the vertex $I$ is admissible
and zero otherwise (see \cite{LW} and \cite{We} for detail).

\vspace{10mm}
\begin{figure}[h]

\includegraphics[height=30mm,width=36mm]{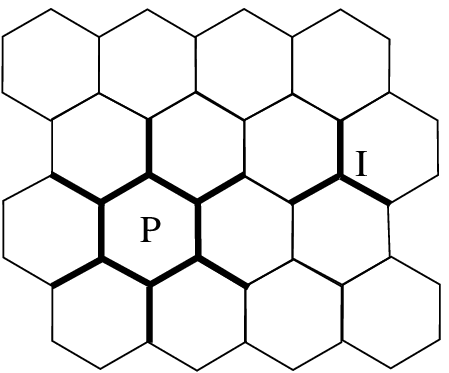}
\hspace{10mm}
\includegraphics[height=30mm,width=30mm]{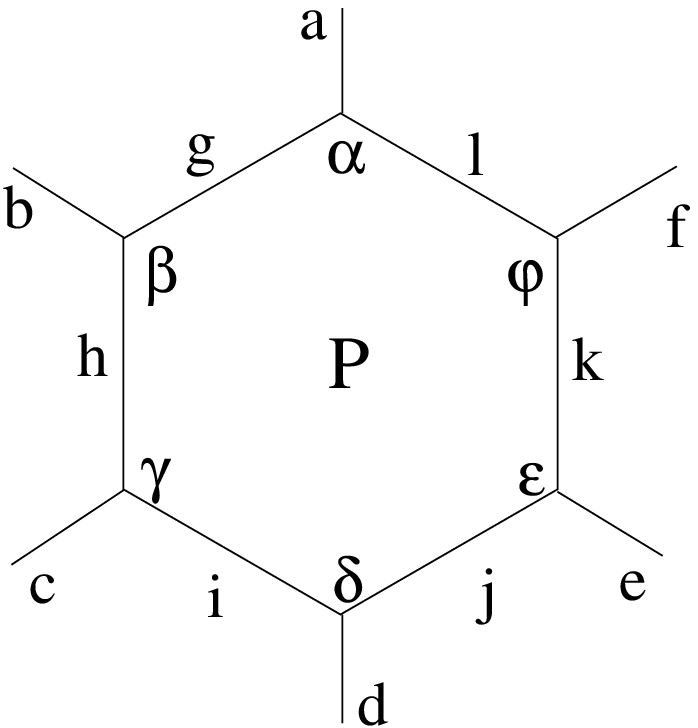}
$\begin{array}{c}\stackrel{B^s_P}{\longmapsto} \\ \\ \\ \\ \\
\\\end{array}$
\includegraphics[height=30mm,width=30mm]{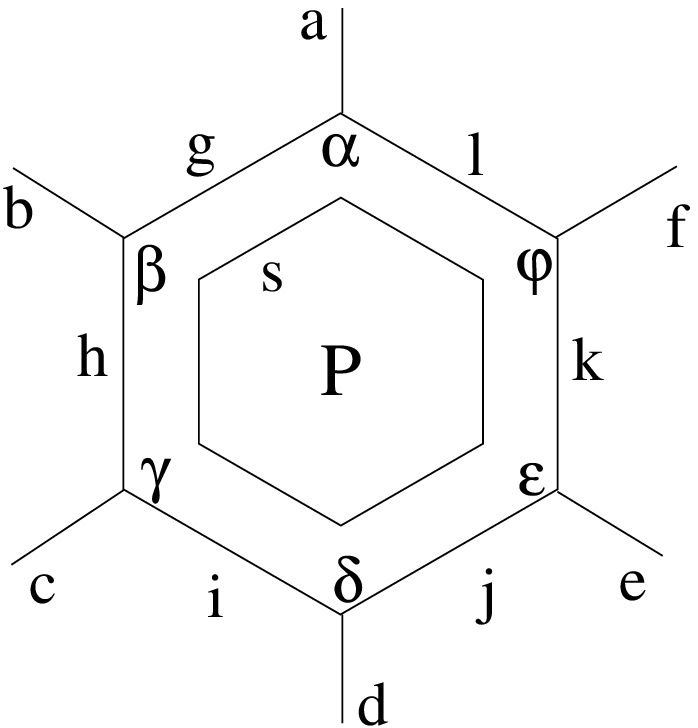}

\caption{Honeycomb lattice model and operator $B^s_P$}
\label{fig:hamiltonian}
\end{figure}

\begin{tabular}{cc}
\end{tabular}

\subsection{$B_P$ is a projector}\label{proj}

The following computation shows that $\left(B_P\right)^2=B_P$:

\vspace{3mm} $\begin{array}{c}(B_P)^2  \\ \\ \end{array}$
\includegraphics[height=12mm,width=15mm]{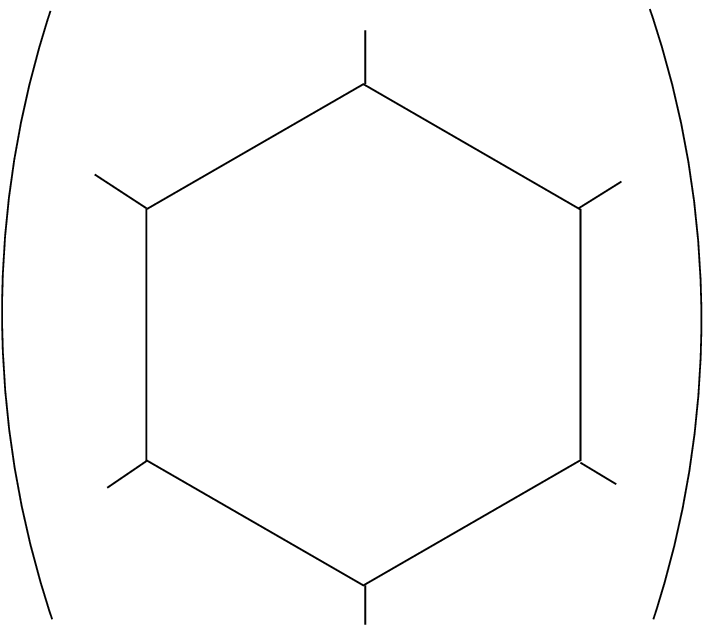}
$\begin{array}{c}= \frac{1}{D^4}\sum_{s,t} st  \\ \\ \end{array}$
\includegraphics[height=12mm,width=12mm]{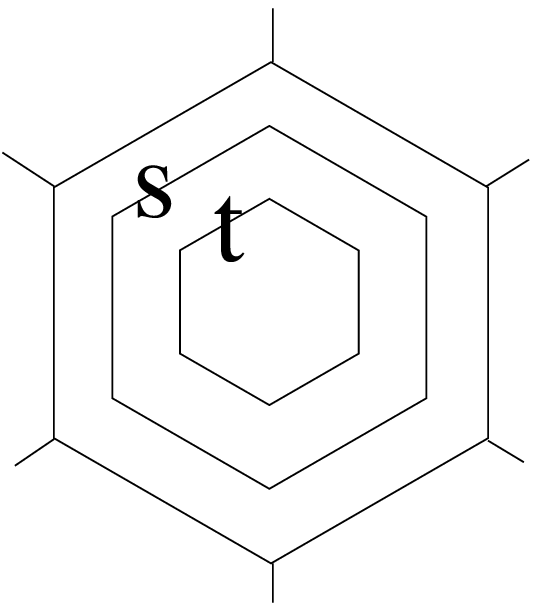}
$\begin{array}{c}= \frac{1}{D^2}\sum_{f} f \\ \\
\end{array}$
\includegraphics[height=12mm,width=12mm]{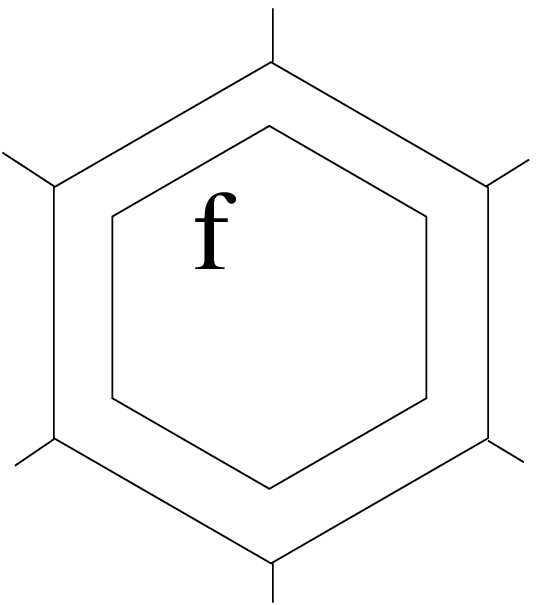}
$\begin{array}{c}= B_P  \\ \\ \end{array}$
\includegraphics[height=12mm,width=15mm]{diagrams/winfig131.eps}

\noindent where the second equality comes from

\vspace{5mm}
\includegraphics[height=12mm,width=10mm]{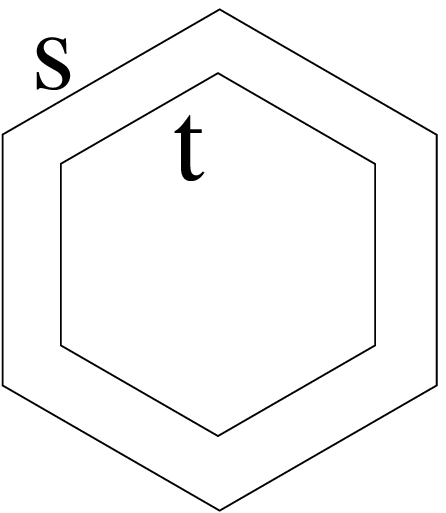}
$\begin{array}{c}=\sum_{f,\eta} \frac{\sqrt{f}}{\sqrt{st}} \\ \\
\end{array}$
\includegraphics[height=12mm,width=12mm]{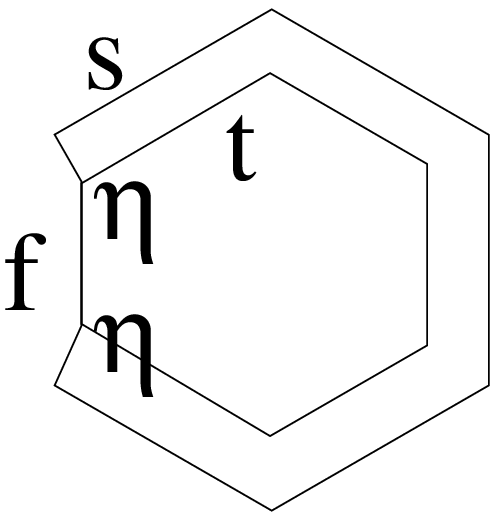}
$\begin{array}{c}=\sum_{f,\eta} \frac{\sqrt{f}}{\sqrt{st}} \\ \\
\end{array}$
\includegraphics[height=12mm,width=13mm]{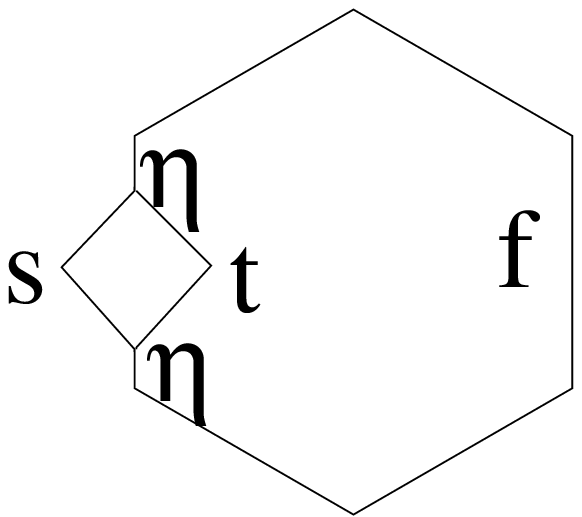}
$\begin{array}{c}=\sum_{f} N^{f}_{st} \\ \\
\end{array}$
\includegraphics[height=12mm,width=10mm]{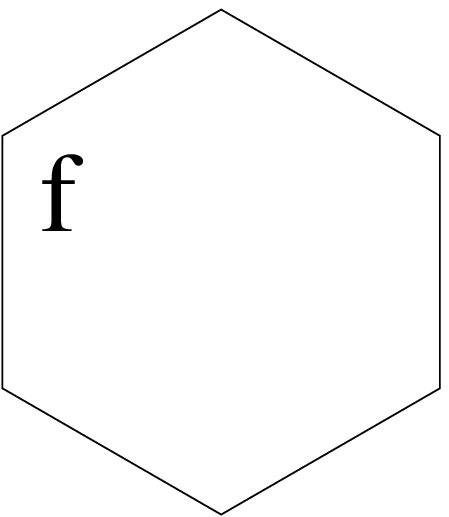},

and $\sum_{s,t}st \sum_{f}N^{f}_{st}$ $=\sum_f \sum_s s
\left(\sum_t tN^{t}_{sf}\right)$ $=\sum_f \sum_s s^2 f$ $=\sum_f
D^2 f$.

\subsection{Hamiltonian $H$ is hermitian}\label{herm}

Let

\vspace{5mm}\hspace{30mm}\includegraphics[height=25mm,width=25mm]{diagrams/winfig113.eps}
$\:\:\: \begin{array}{c}\stackrel{B_P}{\longmapsto} \:\:\: \sum C  \\
\\ \\ \\ \\ \\
\end{array}$
\includegraphics[height=25mm,width=25mm]{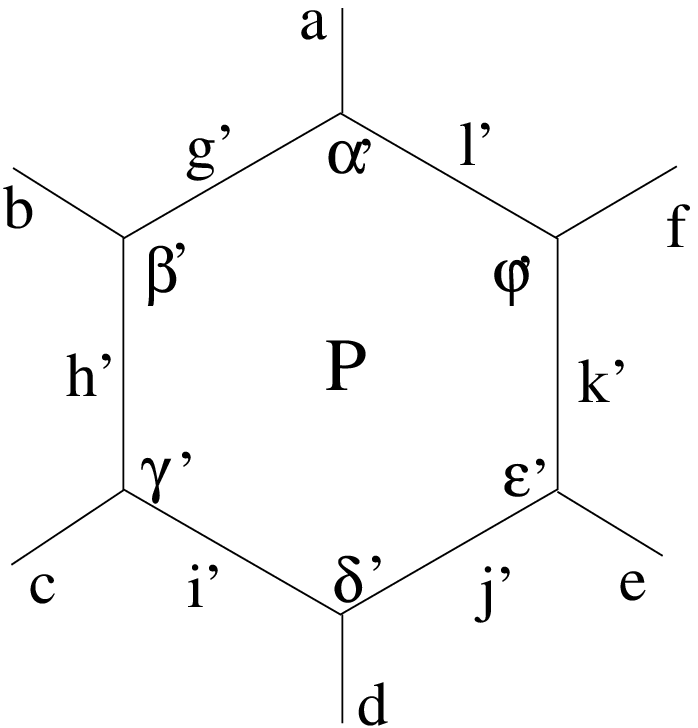}

\noindent where the sum is over all possible
$g',h',i',j',k',l',\alpha',\beta',\gamma',\delta',\varepsilon',\varphi'$
and the coefficient $C$ is a function on 30 labels in both
diagrams. Note that if any label other than the 18 ones around the
hexagonal face is different in both diagrams, then the coefficient
is equal to zero. Once we fix two states shown as above, let us
call them $S$ and $S'$ which are identical outside the diagram.
Then $C=C(S,S')$.

It is sufficient to show that the operator $B_P$ is hermitian,
that is, $C(S,S')$ and $C(S',S)$ are complex conjugate to each
other.

We claim that
$C(S,S')=\frac{1}{D^2}\frac{1}{\sqrt{abcdef}}\chi(S,S')$ where
$\chi(S,S')\in \mathbb{C}$ is defined by the trace of a morphism
in $\Hom_{\mathcal{C}}(d,d)$ as follows:

\vspace{3mm}\hspace{20mm}$\begin{array}{c}\chi(S,S'):=  \\
\\ \\ \\ \\ \\
\end{array}$\includegraphics[height=45mm,width=39mm]{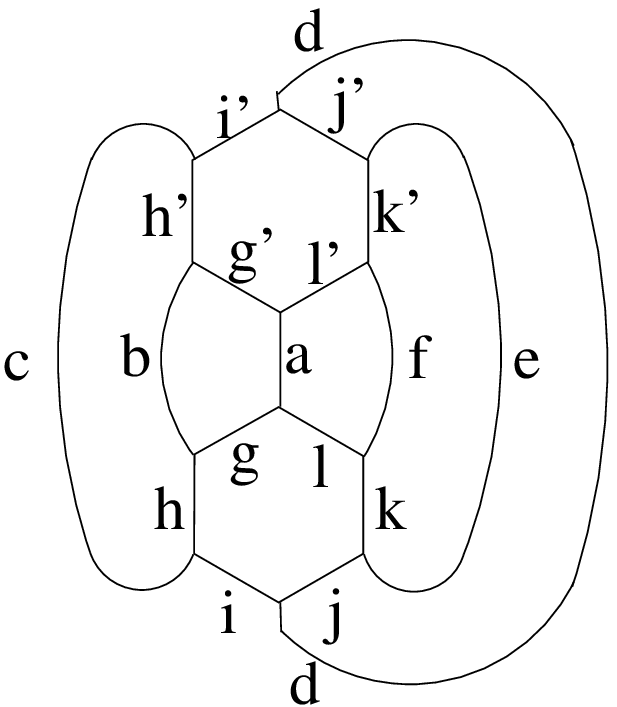}

\noindent We have omitted the labels for trivalent vertices which
are the same as in the states $S$ or $S'$ above determined by the
three edges on each vertex. Note that the hermitian property for
the Hamiltonian is obtained easily from this claim by the mirror
conjugate symmetry since $\chi(S',S)$ is the mirror image of
$\chi(S,S')$. The proof of the claim is done by picture calculus
as below. In this picture calculus, each step indicated by an
arrow has to be a sum of possibly many states, but for the
computation of $C(S,S')$ for fixed states $S$ and $S'$, the
diagram following the arrow is the only state contributing to the
computation, and every other state in the summation is irrelevant.

\vspace{5mm}
\includegraphics[height=25mm,width=25mm]{diagrams/winfig113.eps}
$\begin{array}{c}\:\: \stackrel{(1)}{\rightarrow}\:\:  \\ \\ \\ \\
\end{array}$
\includegraphics[height=25mm,width=25mm]{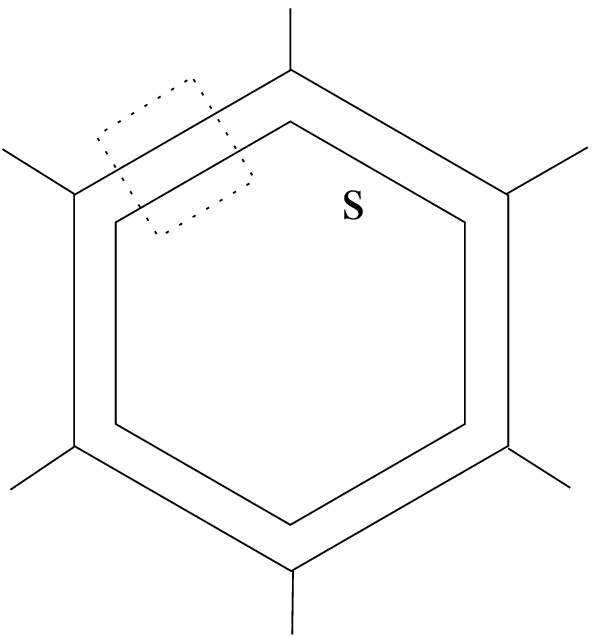}
$\begin{array}{c}\:\: \stackrel{(2)}{\rightarrow}\:\:  \\ \\ \\ \\
\end{array}$
\includegraphics[height=25mm,width=25mm]{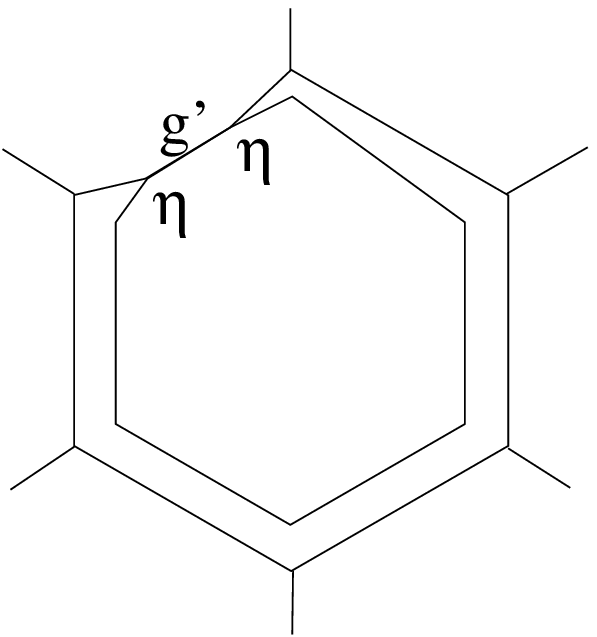}
$\begin{array}{c}\:\: =\:\:  \\ \\ \\
\\\end{array}$
\includegraphics[height=25mm,width=25mm]{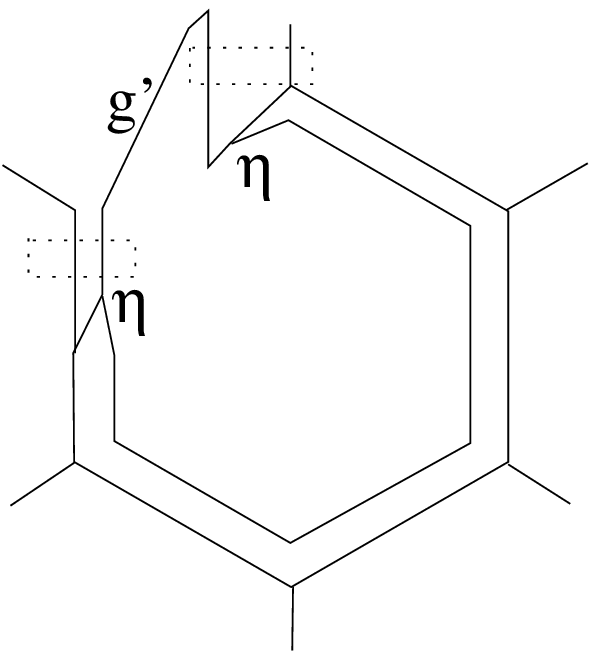}
$\begin{array}{c}\:\: \stackrel{(3)}{\rightarrow}\:\:  \\ \\ \\
\\\end{array}$
\includegraphics[height=25mm,width=25mm]{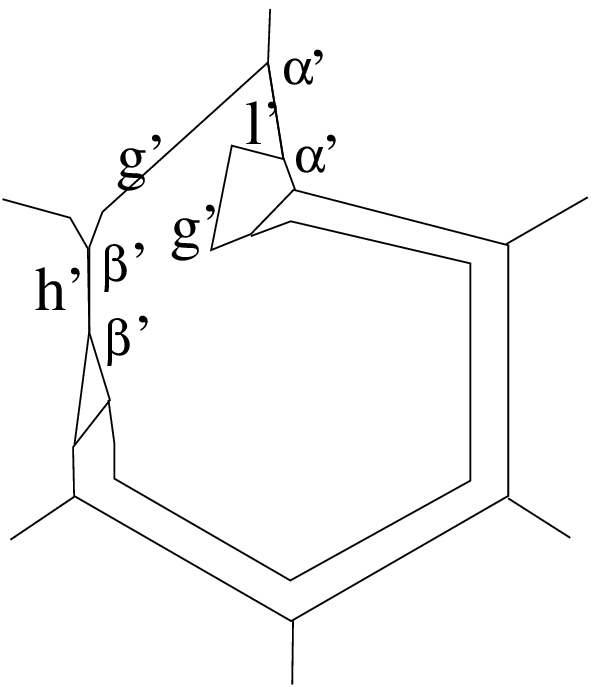}
$\begin{array}{c}\:\: =\:\:  \\ \\ \\
\\\end{array}$
\includegraphics[height=25mm,width=25mm]{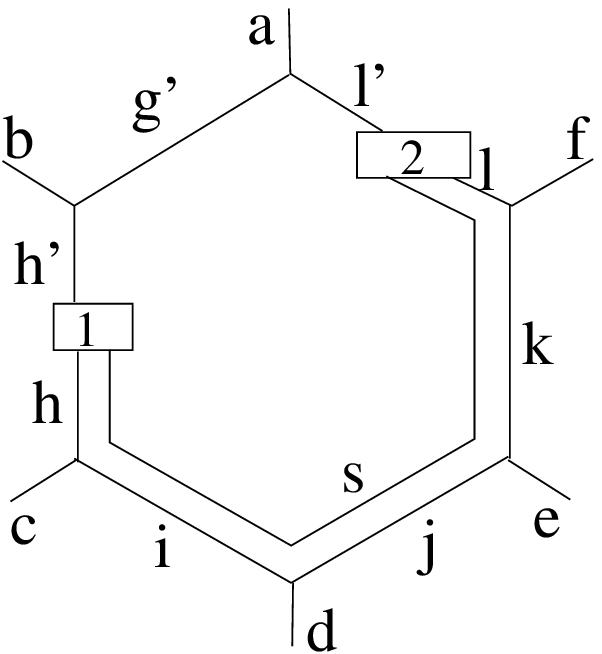}
$\begin{array}{c}\:\: =\:\:  \\ \\ \\
\\\end{array}$
\includegraphics[height=25mm,width=25mm]{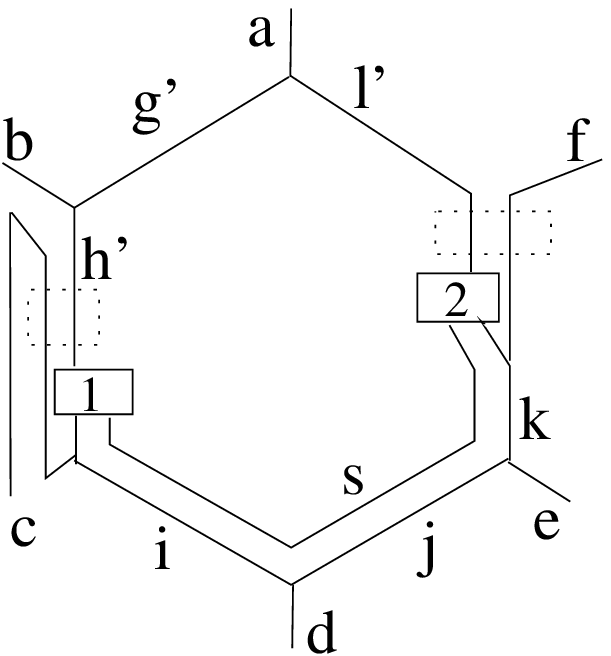}
$\begin{array}{c}\:\: \stackrel{(4)}{\rightarrow}\:\:  \\ \\ \\
\\\end{array}$
\includegraphics[height=25mm,width=25mm]{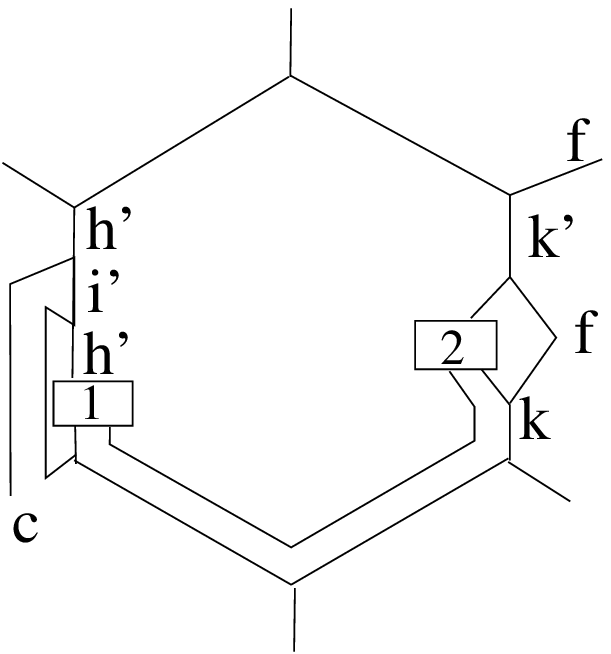}
$\begin{array}{c}\:\: =\:\:  \\ \\ \\
\\\end{array}$
\includegraphics[height=25mm,width=25mm]{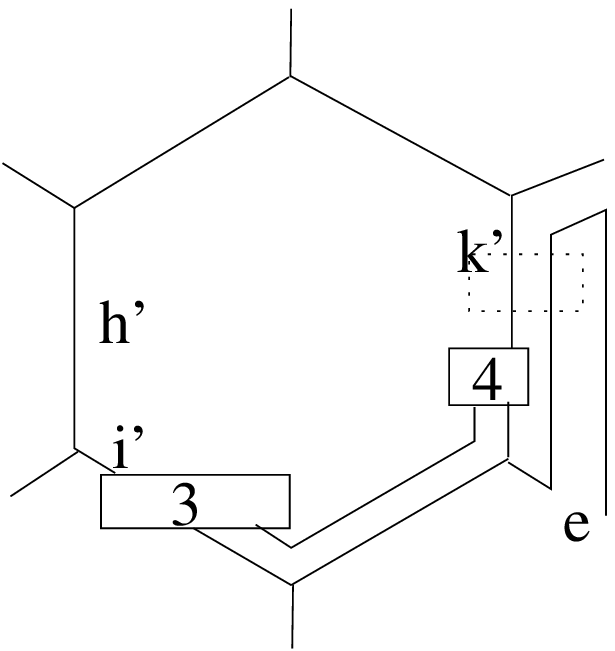}
$\begin{array}{c}\:\: \stackrel{(5)}{\rightarrow}\:\:  \\ \\ \\
\\\end{array}$
\includegraphics[height=25mm,width=25mm]{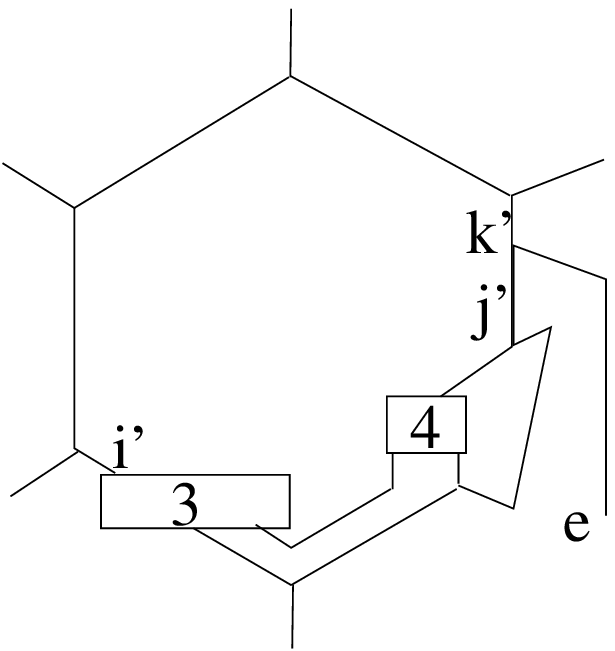}
$\begin{array}{c}\:\: =\:\:  \\ \\ \\
\\\end{array}$
\includegraphics[height=25mm,width=25mm]{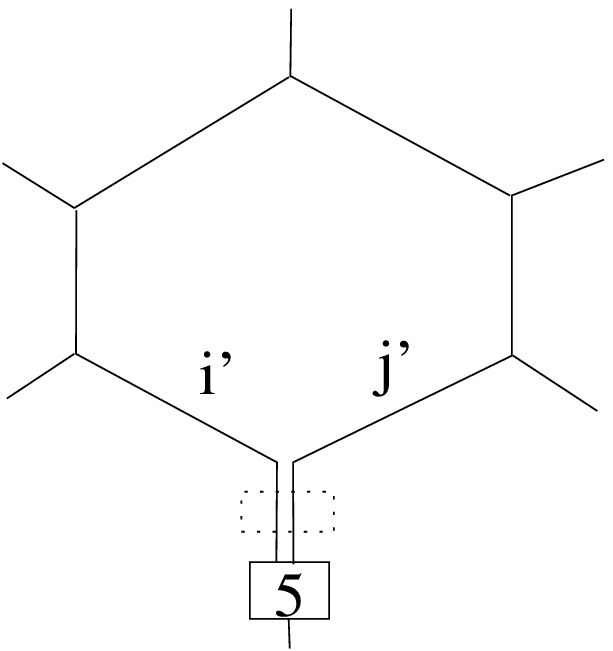}
$\begin{array}{c}\:\: \stackrel{(6)}{\rightarrow}\:\:  \\ \\ \\
\\\end{array}$
\includegraphics[height=25mm,width=25mm]{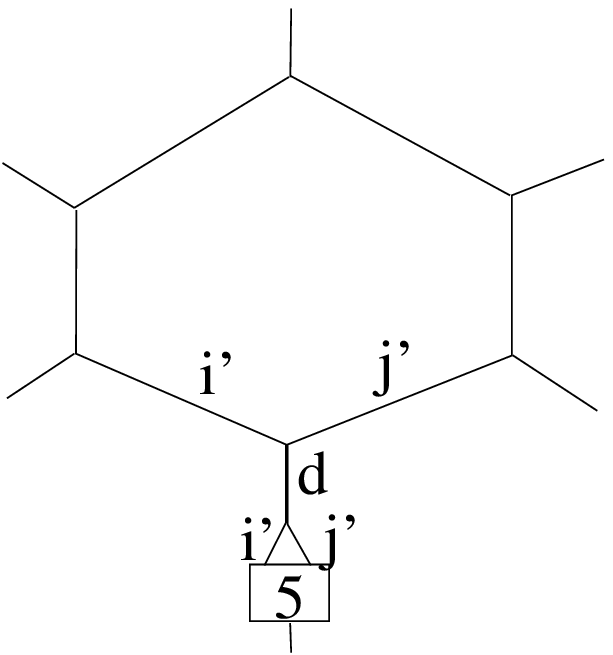}
$\begin{array}{c}\:\: =\:\:  \\ \\ \\
\\\end{array}$
\includegraphics[height=25mm,width=25mm]{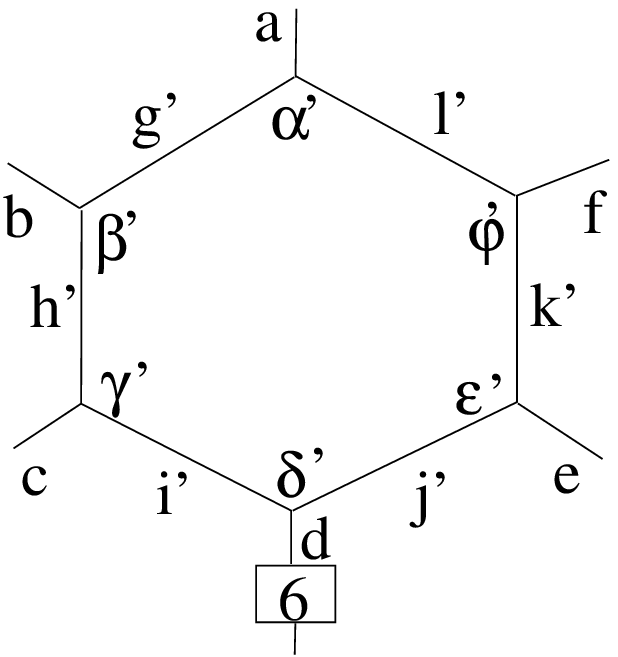}

\noindent where for each step the contributions are as follows:
$(1)=\sum_{s}\frac{s}{D^2}$ applying the operator $B_P$,
$(2)=\sum_{\eta}\frac{\sqrt{g'}}{\sqrt{gs}}$ using (\ref{theta2}),
$(3)=\frac{\sqrt{h'}}{\sqrt{bg'}}\frac{\sqrt{l'}}{\sqrt{g'a}}$
using (\ref{theta2}) and (\ref{theta3}),
$(4)=\frac{\sqrt{i'}}{\sqrt{ch'}}\frac{\sqrt{k'}}{\sqrt{l'f}}$
using (\ref{theta3}), $(5)=\frac{\sqrt{j'}}{\sqrt{k'e}}$ using
(\ref{theta3}), and finally $(6)=\frac{\sqrt{d}}{\sqrt{i'j'}}$
using (\ref{theta2}). So the overall coefficient is
$\left(\frac{\sqrt{d}}{g'\sqrt{abcef}}\frac{1}{D^2}\right)
\frac{\sqrt{g'}}{\sqrt{g}}\sum_{s,\eta}\sqrt{s}$ along with the
diagram in box 6. The following is a deformation of the diagram in
box 6 with scalar
$\frac{\sqrt{g'}}{\sqrt{g}}\sum_{s,\eta}\sqrt{s}$ which completes
the proof.

$\begin{array}{c}\frac{\sqrt{g'}}{\sqrt{g}}\sum_{s,\eta} \sqrt{s} \\ \\
\\
\\\end{array}$
$\begin{array}{c}
\text{\includegraphics[height=15mm,width=7mm]{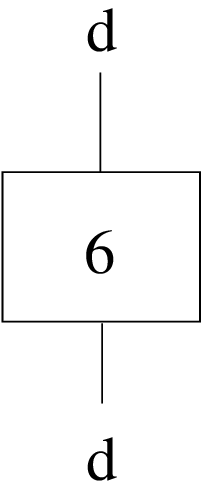}}
\:\:  \\ \\ \\ \\ \end{array}$ $\begin{array}{c} = \frac{\sqrt{g'}}{\sqrt{g}}\sum_{s,\eta} \sqrt{s} \\ \\
\\
\\\end{array}$
\includegraphics[height=30mm,width=33mm]{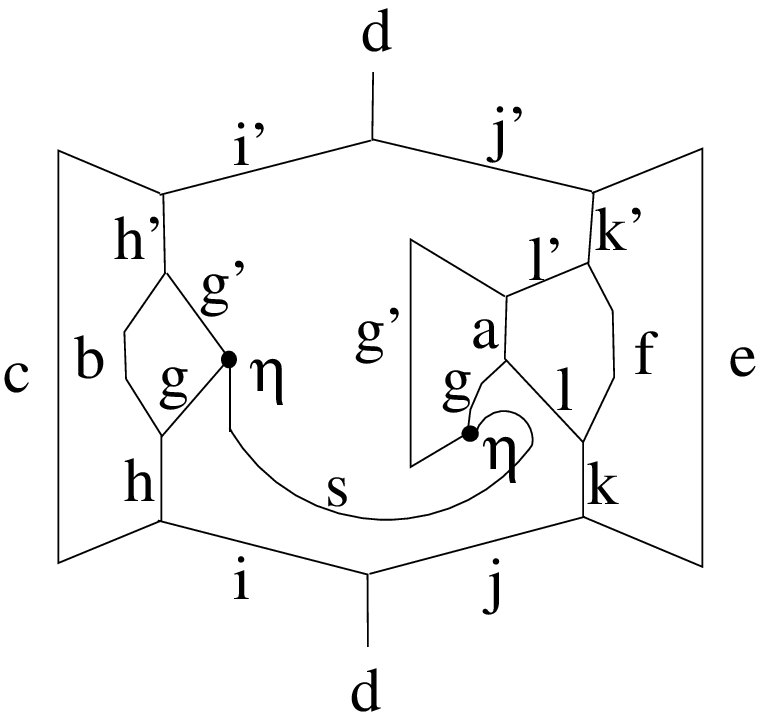}
$\begin{array}{c} = g'  \\ \\
\\
\\\end{array}$
\includegraphics[height=30mm,width=33mm]{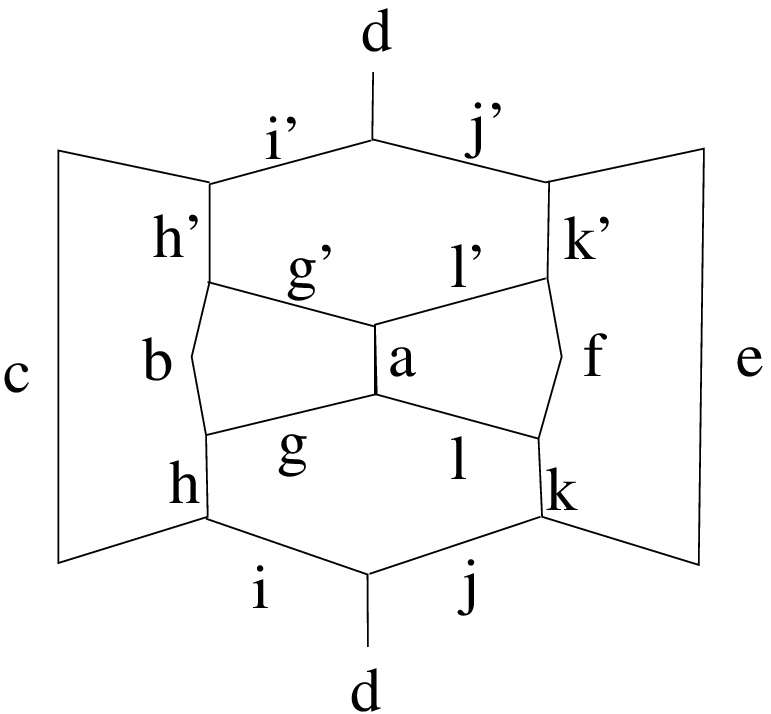}
$\begin{array}{c} = g'\frac{1}{d}\chi(S,S')\id_d  \\ \\
\\
\\\end{array}$

\noindent where the second equality comes from Lemma \ref{can}.

\section{Appendix}

The category $\mathcal{E}$ is a spherical category with three
simple objects, $\{ \textbf{1}, x, y \}$, and fusion rules:

$$x\otimes y=y\otimes x=x, x\otimes x=\textbf{1}\oplus 2x \oplus y, y\otimes y=\textbf{1}$$

The following is the list of $F$-matrices obtained by transposing
the ones in \cite{HH} since we are using left multiplication
convention in this paper.

\vspace{3mm}

$F^{y}_{yyy}=$ $F^{x}_{xyy}=$ $F^{x}_{yyx}=$
$F^{\textbf{1}}_{xyx}=$ $F^{\textbf{1}}_{xxy}=$ $F^{y}_{xxy}=$
$F^{\textbf{1}}_{yxx}=$ $F^{y}_{yxx}=1$,

$F^{y}_{xyx}=$ $F^{x}_{yxy}= -1$,

$F^{x}_{xyx}= \left[\begin{array}{cc}1&0\\0&-1\end{array}
\right]$, $F^{x}_{xxy}=
\left[\begin{array}{cc}0&-i\\i&0\end{array} \right]$,
$F^{x}_{yxx}= \left[\begin{array}{cc}0&1\\1&0\end{array} \right]$,
$F^{\textbf{1}}_{xxx}= \frac{1}{\sqrt{2}} e^{7 \pi i
 /12} \left[\begin{array}{cc}1&1\\i&-i\end{array}
\right]$, $F^{y}_{xxx}= \frac{1}{\sqrt{2}} e^{7 \pi i
 /12} \left[\begin{array}{cc}i&-i\\1&1\end{array}
\right]$,

$F^{x}_{xxx}=$ $\left[ \begin{smallmatrix}
\frac{-1+\sqrt{3}}{2}&\frac{-1+\sqrt{3}}{2}&1&1&1&-1\\
\frac{-1+\sqrt{3}}{2}&\frac{1-\sqrt{3}}{2}&1&1&-1&1\\
\frac{1- \sqrt{3}}{4} e^{\pi i /6}&\frac{1- \sqrt{3}}{4} e^{\pi i /6}&-\frac{1}{2} (e^{\pi i /6} -1)&\frac{1}{2} e^{\pi i /3}&-\frac{1}{2} (e^{\pi i /6} -1)&-\frac{1}{2} e^{\pi i /3}\\
\frac{1- \sqrt{3}}{4} e^{2\pi i /3}&\frac{1- \sqrt{3}}{4} e^{2\pi i /3}&\frac{1}{2} e^{5\pi i /6}&\frac{1}{2} (e^{-\pi i /3} + i)&\frac{1}{2} e^{5\pi i /6}&-\frac{1}{2} (e^{-\pi i /3} + i)\\
\frac{1- \sqrt{3}}{4} e^{2\pi i /3}&-\frac{1- \sqrt{3}}{4} e^{2\pi
i /3}&\frac{1}{2} (e^{-\pi i /3} + i)&\frac{1}{2} e^{5\pi i
/6}&-\frac{1}{2} (e^{-\pi i /3} +
i)&\frac{1}{2} e^{5\pi i /6}\\
\frac{1- \sqrt{3}}{4} e^{\pi i /6}&-\frac{1- \sqrt{3}}{4} e^{\pi i
/6}&\frac{1}{2} e^{\pi i /3}&-\frac{1}{2} (e^{\pi i /6}
-1)&-\frac{1}{2} e^{\pi i /3}&-\frac{1}{2} (e^{\pi i /6} -1)
\end{smallmatrix}
\right]$

$\left(F^{x}_{xxx}\right)^{-1}=$ $\left[ \begin{smallmatrix}
\frac{1}{1+\sqrt{3}}&\frac{1}{1+\sqrt{3}}&e^{5\pi i/6}&e^{\pi i/3}&e^{\pi i/3}&e^{5\pi i/6}\\
\frac{1}{1+\sqrt{3}}&-\frac{1}{1+\sqrt{3}}&e^{5\pi i/6}&e^{\pi i/3}&e^{-2\pi i/3}&e^{-\pi i/6}\\
\frac{1}{2(1+\sqrt{3})}&\frac{1}{2(1+\sqrt{3})}&\frac{1}{\sqrt{2}(1+\sqrt{3})}e^{5\pi i/12}&\frac{1}{2}e^{-5\pi i/6}&\frac{1}{\sqrt{2}(1+\sqrt{3})}e^{-\pi i/12}&\frac{1}{2}e^{-\pi i/3}\\
\frac{1}{2(1+\sqrt{3})}&\frac{1}{2(1+\sqrt{3})}&\frac{1}{2}e^{-\pi i/3}&\frac{1}{\sqrt{2}(1+\sqrt{3})}e^{-\pi i/12}&\frac{1}{2}e^{-5\pi i/6}&\frac{1}{\sqrt{2}(1+\sqrt{3})}e^{5\pi i/12}\\
\frac{1}{2(1+\sqrt{3})}&-\frac{1}{2(1+\sqrt{3})}&\frac{1}{\sqrt{2}(1+\sqrt{3})}e^{5\pi i/12}&\frac{1}{2}e^{-5\pi i/6}&\frac{1}{\sqrt{2}(1+\sqrt{3})}e^{11\pi i/12}&\frac{1}{2}e^{2\pi i/3}\\
-\frac{1}{2(1+\sqrt{3})}&\frac{1}{2(1+\sqrt{3})}&\frac{1}{2}e^{2\pi
i/3}&\frac{1}{\sqrt{2}(1+\sqrt{3})}e^{11\pi
i/12}&\frac{1}{2}e^{-5\pi
i/6}&\frac{1}{\sqrt{2}(1+\sqrt{3})}e^{5\pi i/12}
\end{smallmatrix}
\right]$

\end{document}